\documentclass{conm-p-l}

\copyrightinfo{2004}{}

\newtheorem{theorem}{Theorem}[section]
\newtheorem{lemma}[theorem]{Lemma}
\newtheorem{corollary}[theorem]{Corollary}
\theoremstyle{definition}

\theoremstyle{remark}
\newtheorem{remark}[theorem]{Remark}

\numberwithin{equation}{section}

\newcommand{\D}{{\mathcal D}}
\newcommand{\tla}{\tilde{\lambda}}

\newcommand{\tu}{\tilde{u}}

\newcommand{\hk}{\hat{k}}
\newcommand{\Zo}{\mathbb{Z}/\{0\}}
\newcommand{\ZZo}{\mathbb{Z}^2/\{0\}}
\newcommand{\Z}{\mathbb{Z}}
\newcommand{\R}{\mathbb{R}}
\newcommand{\C}{\mathbb{C}}

\newcommand{\e}{\epsilon}

\newcommand{\ta}{\tilde{a}}

\newcommand{\ga}{\gamma}
\newcommand{\Ga}{\Gamma}

\newcommand{\dl}{\delta}
\newcommand{\Dl}{\Delta}

\newcommand{\ra}{\rightarrow}
\newcommand{\al}{\alpha}
\newcommand{\be}{\beta}

\newcommand{\pa}{\partial}

\newcommand{\La}{\Lambda}

\newcommand{\la}{\lambda}

\newcommand{\nid}{\noindent}

\newcommand{\tom}{\tilde{\omega}}
\newcommand{\om}{\omega}
\newcommand{\Om}{\Omega}
\newcommand{\bOm}{\bar{\Omega}}
\newcommand{\na}{\nabla}

\newcommand{\non}{\nonumber}

\begin{document}

\title[Invariant Manifolds for Navier-Stokes Equations]
{Invariant Manifolds and Their Zero-Viscosity Limits for 
Navier-Stokes Equations}

\author{Y. Charles Li}
\address{Department of Mathematics, University of Missouri, 
Columbia, MO 65211}
\curraddr{}
\email{cli@math.missouri.edu}
\thanks{}


\subjclass{Primary 35, 76, 37; Secondary 34}
\date{}

\dedicatory{}


\begin{abstract}
First we prove a general spectral theorem for the linear Navier-Stokes 
(NS) operator in both 2D and 3D. The spectral theorem says that the 
spectrum consists of only eigenvalues which lie in a parabolic region, 
and the eigenfunctions (and higher order eigenfunctions) form a complete 
basis in $H^\ell$ ($\ell = 0,1,2, \cdots$). Then we prove the existence 
of invariant manifolds. We are also interested in a more challenging 
problem, i.e. studying the zero-viscosity limits ($\nu \ra 0^+$) of the 
invariant manifolds. Under an assumption, we can show that the sizes of 
the unstable manifold and the center-stable manifold of a steady state 
are $O(\sqrt{\nu})$, while the sizes of the stable manifold, the center 
manifold, and the center-unstable manifold are $O(\nu)$, as $\nu \ra 0^+$. 
Finally, we study three examples. The first example is defined on a 
rectangular periodic domain, and has only one unstable eigenvalue which 
is real. A complete estimate on this eigenvalue is obtained. Existence 
of an 1D unstable manifold and a codim 1 stable manifold is proved without
any assumption. For the other two examples, partial estimates on the 
eigenvalues are obtained.
\end{abstract}

\maketitle

\tableofcontents








\section{Introduction}

Navier-Stokes equations define an infinite dimensional dynamical 
system. Turbulence can be represented through a chaotic solution 
in the infinite dimensional phase spaces. For example, one can take 
the Sobolev spaces $H^\ell$ ($\ell = 0,1,2, \cdots$) to be the 
phase spaces. A basic dynamical system question is: Does the 
phase space foliate under the NS flow ? That is, are there invariant 
manifolds under the NS flow ? The answer is yes. It seems that 
under the Euler flow, which can be obtained from the NS flow by 
setting viscosity to zero, the phase space does not foliate. To begin 
a dynamical system study, one chooses a common steady state of the NS and its 
corresponding Euler equation, defined on a periodic spatial domain. At 
this steady state, 
the spectra of the linear NS and the linear Euler are dramatically 
different. The linear NS has the form: $\nu \Dl +$ relatively compact 
terms, therefore, it has only a point spectrum via 
Weyl's essential spectrum theorem. On the other hand, the linear Euler 
has a nontrivial essential spectrum. For 2D linear Euler, if the steady 
state as a vector field has a saddle, then the essential spectrum in the 
Sobolev space $H^\ell$ ($\ell =0,1,2, \cdots$) is a vertical band symmetric 
with respect to the imaginary axis \cite{SL03}. The width of the band 
($=\ell \La $ for some $\La > 0$) increases with $\ell$, representing 
the fact of cascade and inverse cascade. A typical example of such steady 
states is the cat's eye. If the steady state has no saddle, but has at 
least two centers, then the essential spectrum in the Sobolev space 
$H^\ell$ is the imaginary axis \cite{SL03}. A typical example of such 
steady states is the shear \cite{MS61} \cite{Li00}. For the shear, a 
sharp upper bound 
on the number of eigenvalues was also obtained \cite{LLS04}. For 3D linear 
Euler, the spectrum is an open problem. A nice example of the steady states 
is the ABC flow vector field.

An interesting open problem is to understand the connection between the 
spectra of linear NS and linear Euler as $\nu \ra 0^+$. Besides the 
progresses mentioned above in understanding the spectra of linear Euler, 
progresses were also made in better understanding of the spectra of linear 
NS, although they were in a different setting. For both 2D and 3D NS posed 
in $L^2$ (in terms of weak solutions) on a bounded domain with non-slip 
boundary condition, Prodi \cite{Pro62} proved that the eigenvalues of the 
linear NS at a steady state lie inside a parabolic region 
$\la_r < a-b \la_i^2$ where $a>0$ and $b>0$ are some constants, and 
$\la = \la_r + i \la_i$ is the spectral parameter. Moreover, David 
Sattinger \cite{Sat70} proved that the linear NS has infinitely many 
eigenvalues lying in the parabolic region, each of finite multiplicity, 
which can cluster only at infinity, and the corresponding eigenfunctions 
(and higher order eigenfunctions) form a complete basis in $L^2$. Similar 
results are proved in the current article for a periodic domain and mild 
solutions in $H^\ell$ ($\ell =0,1,2, \cdots$). The proof involves some 
deep result from Hilbert-Schmidt operators. Together with some deep result 
from sectorial operator, we prove the existence of invariant manifolds.

In terms of well-posedness of Navier-Stokes equations (NS), 
zero-viscosity limit problem has been extensively studied 
\cite{Kat75} \cite{Wu96}. Here we are also interested in the zero-viscosity 
limit ($\nu \ra 0^+$) problem for invariant manifolds.
It is proved in this article that 
under an assumption, their sizes are of order $O(\sqrt{\nu})$ 
or $O(\nu )$ as $\nu \ra 0^+$. Existence of invariant manifolds 
for Euler equations is still open. I tend to believe that invariant manifolds
for Euler equations do not exist, at least not in the usual sense.
An intuitive argument is as follows: Kato's technique \cite{Kat75} for 
proving the existence of a mild solution to Euler equations amounts 
to writing the solution as $u(t) = e^{-\int_0^t u(\tau )\cdot \na d \tau }
u(0)$, where $u$ is the velocity, and the Leray projection and the external 
force are not included, which do not affect the argument. Notice that 
the $u(\tau )$ in the exponent depends upon $u(0)$, therefore, one 
should only expect the evolution operator of Euler equations to be 
at best $C^0$ in $u(0)$, just as in $t$. On the other hand, the existence 
of an invariant manifold is really an at least $C^1$ in $u(0)$ phenomenon.
One can obtain expressions of invariant manifolds around a shear for a 
Galerkin truncation of the Euler equation \cite{Li03}, which exhibit a lip 
structure. 

For concrete examples, e.g. single mode steady states, a continued fraction 
technique can be designed for calculating the eigenvalues \cite{MS61}. 
Detailed calculations were conducted by Vincent Liu \cite{Liu92} \cite{Liu95}.
Liu obtained rather detailed information on unstable eigenvalues. In some 
case, Liu showed that as $\nu \ra 0^+$, the unstable eigenvalue does 
converge to that of linear Euler. In this article, besides the example of Liu, 
we will study two more examples. For the first example, a complete eigenvalue
information is obtained, and the existence of unstable and stable manifolds 
is proved without any assumption. For the second example, only partial 
information on the eigenvalues is obtained, as in the example of Liu.

The article is organized as follows: In section 2, we will present the 
formulation of the problems. Section 3 is on invariant manifolds. Examples are 
studied in Section 4. 

\section{Formulation of the Problems}

We will study the following form of 2D Navier-Stokes equation (NS)
\begin{equation}
\pa_t \Om + \{ \Psi, \Om \} = \nu [\Dl \Om + f(x)] \ ,
\label{2DNS}
\end{equation}
where $\Om$ is the vorticity which is a real scalar-valued function
of three variables $t$ and $x=(x_1, x_2)$, the bracket $\{\ ,\ \}$ 
is defined as
\[
\{ f, g\} = (\pa_{x_1} f) (\pa_{x_2}g) - (\pa_{x_2} f) (\pa_{x_1} g) \ ,
\]
where $\Psi$ is the stream function given by,
\[
u_1=- \pa_{x_2}\Psi \ ,\ \ \ u_2=\pa_{x_1} \Psi \ ,
\]
the relation between vorticity $\Om$ and stream 
function $\Psi$ is,
\[
\Om =\pa_{x_1} u_2 - \pa_{x_2} u_1 =\Dl \Psi \ ,
\]
and $\nu$ is the viscosity, $\Dl$ is the 2D Laplacian , and 
$f(x)$ is the body force. We will pose periodic boundary condition
\[
\Om (t, x_1 +2\pi , x_2) = \Om (t, x_1 , x_2) = \Om (t, x_1, x_2 +2\pi ),
\]
i.e. the 2D NS is defined on the 2-torus $\mathbb{T}^2$. We require that 
$\Psi$ and $f$ have mean zero
\[
\int_{\mathbb{T}^2} \Psi dx = \int_{\mathbb{T}^2} f dx = 0\ .
\]
Of course $\Om$ always has mean zero. In this case, $\Psi = \Dl^{-1} \Om $.

Setting $\nu = 0$ in the 2D NS (\ref{2DNS}), one gets the corresponding 
2D Euler equation. We will study the common steady states of 2D NS and 
2D Euler given by 
\begin{equation}
\{ \Psi, \Om \} = 0 \ , \quad  \Dl \Om + f(x) = 0 \ . 
\label{2DSS}
\end{equation}
Our three prime examples are: (1). The shear $\Om = \cos x_2$ defined 
on the rectangular periodic domain $[0,2\pi/\al ]\times [0, 2\pi ]$ 
where $1/2 < \al < 1$. (2). The shear $\Om = \cos (x_1+x_2)$ where 
$f(x)=2 \cos (x_1+x_2)$. (3). The cat's eye $\Om = \cos x_1 +\al 
\cos x_2$ where $\al$ is a constant and $f(x)=\cos x_1 +\al \cos x_2$.

We will study the following form of 3D Navier-Stokes equation
\begin{equation}
\pa_t \Om + (u \cdot \na) \Om - (\Om \cdot \na) u = \nu [\Dl \Om + f(x)] \ ,
\label{3DNS}
\end{equation}
where $u = (u_1, u_2, u_3)$ is the velocity, $\Om = (\Om_1, \Om_2, \Om_3)$
is the vorticity, $\na = (\pa_{x_1}, \pa_{x_2}, \pa_{x_3})$, 
$\Om = \na \times u$, $\na \cdot u = 0$, $\nu$ is the viscosity, $\Dl$ is 
the 3D Laplacian, and $f(x) = (f(x_1), f(x_2), f(x_3))$ is the body force. 
We also pose periodic boundary condition of period ($2\pi, 2\pi, 2\pi$), i.e.
the 3D NS is defined on the 3-torus $\mathbb{T}^3$. We require that $u$, 
$f$ and $\Om$ all have mean zero. In this case, $u$ can be uniquely 
determined from $\Om$ by Fourier transform:
\begin{eqnarray*}
\tu_1(k) &=& i |k|^{-2} [ k_2 \tilde{\Om}_3(k) - k_3 \tilde{\Om}_2(k)]\ , \\
\tu_2(k) &=& i |k|^{-2} [ k_3 \tilde{\Om}_1(k) - k_1 \tilde{\Om}_3(k)]\ , \\
\tu_3(k) &=& i |k|^{-2} [ k_1 \tilde{\Om}_2(k) - k_2 \tilde{\Om}_1(k)]\ .
\end{eqnarray*}

Setting $\nu = 0$ in the 3D NS (\ref{3DNS}), one gets the corresponding 
3D Euler equation. We will study the common steady states of 3D NS and 
3D Euler given by 
\begin{equation}
(u \cdot \na) \Om - (\Om \cdot \na) u = 0 \ , \quad 
\Dl \Om + f(x) = 0 \ .
\label{3DSS}
\end{equation}
Our prime example is: The ABC flow
\[
u_1 = A \sin x_3 + C \cos x_2 \ , \quad 
u_2 = B \sin x_1 + A \cos x_3 \ , \quad
u_3 = C \sin x_2 + B \cos x_1 \ . 
\]
\nid
Here $\Om = u = f$.

\section{Invariant Manifolds}

Centering around a $C^\infty$ steady state $\Om_*$, $u_*$, $\Psi_*$ given 
by (\ref{2DSS}) or (\ref{3DSS}), one can rewrite the NS as 
\begin{equation}
\pa_t \Om = L \Om + N(\Om)\ ,
\label{NSr}
\end{equation}
where for 2D 
\begin{eqnarray*}
L\Om &=& \nu \Dl \Om - \{ \Psi_*, \Om \} - \{ \Psi, \Om_* \} \ , \\
N(\Om) &=& -\{ \Psi, \Om \} \ ,
\end{eqnarray*}
and for 3D 
\begin{eqnarray*}
L\Om &=& \nu \Dl \Om - (u_* \cdot \na )\Om - (u \cdot \na )\Om_* \\
     &+& (\Om_* \cdot \na )u +(\Om \cdot \na )u_* \ , \\
N(\Om) &=& - (u \cdot \na )\Om +(\Om \cdot \na )u\ . 
\end{eqnarray*}
\begin{theorem}
The linear Euler operator $L$ is a closed operator. Its spectrum 
consists of infinitely many discrete eigenvalues, each of finite 
multiplicity, which can cluster only at infinity. All the eigenvalues lie
inside the parabolic region
\[
\la_r \leq a-b\la_i^2
\]
where $a>0$ and $b>0$ are some constants, and $\la = \la_r + i \la_i$ is 
the spectral parameter. The corresponding eigenfunctions (and higher order
eigenfunctions) form a complete basis in $H^\ell$ ($\ell = 0,1,2, \cdots$).
\end{theorem}
\begin{proof}
The proof for $L$ being a closed operator is trivial.
We can write $L$ as 
\[
L = \nu \Dl + A
\]
where $A$ is relatively compact with respect to $\nu \Dl$. By Weyl's 
essential spectrum theorem \cite{RS78}, $L$ and $\nu \Dl$ have the same 
essential spectrum. Since the essential spectrum of $\nu \Dl$ is empty, 
so is $L$. Therefore, $L$ has only eigenvalues. It is easy to see that 
for large enough $\la >0$, $(L-\la )^{-1}$ is compact, thus each 
eigenvalue can be of only finite multiplicity and all the eigenvalues 
can cluster only at infinity \cite{BJS64}. Next we show that large 
enough $\la >0$ is not an eigenvalue. Assume the contrary, then there exists
nonzero $\Om \in H^{\ell +2}$ such that
\[
\la \Om -\nu \Dl \Om = A\Om \ .
\]
Applying $(-\Dl)^{\ell /2}$, one gets
\[
\la (-\Dl)^{\ell /2}\Om -\nu \Dl (-\Dl)^{\ell /2}\Om =(-\Dl)^{\ell /2} 
(A\Om )\ .
\]
Multiplying by $(-\Dl)^{\ell /2} \Om$, one obtains
\begin{eqnarray*}
& & \la \left |(-\Dl)^{\ell /2}\Om \right |^2 + \nu \left [(-\Dl)^{\ell /2} \Om
\right ] \cdot \left [ (-\Dl) (-\Dl)^{\ell /2}\Om \right ] \\
&=& \left [(-\Dl)^{\ell /2} \Om \right ] \cdot \left [ (-\Dl)^{\ell /2} (A\Om )
\right ]  \ .
\end{eqnarray*}
Thus
\begin{eqnarray*}
& & \la \int \left |(-\Dl)^{\ell /2}\Om \right |^2 + \nu \int 
\left | (-\Dl)^{ (\ell +1) /2} \Om \right |^2 \\
&=&  \int
\left [(-\Dl)^{\ell /2} \Om \right ] \cdot \left [ (-\Dl)^{\ell /2} (A\Om )
\right ]   \\
&\leq& C \int \left | (-\Dl)^{ \ell /2} \Om \right |
\left | (-\Dl)^{ (\ell +1) /2} \Om \right | \\
&\leq& C \left ( \e^2 \int \left | (-\Dl)^{ (\ell +1) /2} \Om \right |^2 
+\frac{1}{\e^2} \int \left |(-\Dl)^{\ell /2}\Om \right |^2 \right )\ .
\end{eqnarray*}
Let $\e^2 = \frac{\nu}{2C}$ and $\la = \frac{2C^2}{\nu} +1$, we have 
\[
\int \left |(-\Dl)^{\ell /2}\Om \right |^2 + \frac{1}{2} \nu \int 
\left | (-\Dl)^{ (\ell +1) /2} \Om \right |^2 \leq 0 \ ,
\]
which is a contradiction to $\Om$ being an eigenfunction. Thus 
$\la = \frac{2C^2}{\nu} +1$ is in the resolvent set. Next let $\la$ be 
a complex number in the resolvent set. Naturally, we need to consider $\Om$ 
as a complex-valued function. For any $h \in H^\ell$, there is a $\Om 
\in H^{\ell +2}$ such that
\[
-\nu \Dl \Om + \la \Om - A\Om = h \ .
\]
Applying $(-\Dl)^{\ell /2}$, one gets
\[
-\nu \Dl (-\Dl)^{\ell /2}\Om +\la (-\Dl)^{\ell /2}\Om - (-\Dl)^{\ell /2} 
(A\Om ) = (-\Dl)^{\ell /2} h \ .
\]
Multiplying by $(-\Dl)^{\ell /2} \bOm $, one obtains
\begin{eqnarray*}
& & \nu \left [ (-\Dl) (-\Dl)^{\ell /2}\Om \right ] \cdot
\left [(-\Dl)^{\ell /2} \bOm \right ] + \la 
\left [(-\Dl)^{\ell /2} \Om \right ] \cdot  
\left [(-\Dl)^{\ell /2} \bOm \right ] \\
&=& \left [ (-\Dl)^{\ell /2} (A\Om )\right ] \cdot  
\left [(-\Dl)^{\ell /2} \bOm \right ] + 
\left [(-\Dl)^{\ell /2} h \right ] \cdot  
\left [(-\Dl)^{\ell /2} \bOm \right ] \ .
\end{eqnarray*}
Thus
\begin{eqnarray}
& & \nu \int \left | (-\Dl)^{ (\ell +1) /2} \Om \right |^2 +
\la \int \left |(-\Dl)^{\ell /2}\Om \right |^2  \nonumber \\
&=&  \int \left [ (-\Dl)^{\ell /2} (A\Om )\right ] \cdot  
\left [(-\Dl)^{\ell /2} \bOm \right ] + \int 
\left [(-\Dl)^{\ell /2} h \right ] \cdot  
\left [(-\Dl)^{\ell /2} \bOm \right ] \ . \label{SDS1}
\end{eqnarray}
Taking the real part, one gets
\begin{eqnarray*}
\nu \| \Om \|_{\ell +1}^2 +\la_r \| \Om \|^2_\ell &\leq& C 
\left (\e^2 \| \Om \|_{\ell +1}^2 + \frac{1}{\e^2} \| \Om \|^2_\ell \right ) \\
&+&\frac{1}{2}\left (\| h \|^2_\ell + \| \Om \|^2_\ell \right ) \ .
\end{eqnarray*}
Let $\e^2 = \frac{\nu}{2C}$, one has 
\begin{equation}
\frac{1}{2} \nu \| \Om \|_{\ell +1}^2 +\la_r \| \Om \|^2_\ell
\leq \frac{1}{2}\| h \|^2_\ell + C_1 \| \Om \|^2_\ell \ .
\label{SDS2}
\end{equation}
Taking the imaginary part of (\ref{SDS1}), one gets
\[
\la_i \| \Om \|^2_\ell \leq C_2 \| \Om \|_{\ell +1} \| \Om \|_\ell
+ \| h \|_\ell \| \Om \|_\ell \ ,
\]
that is
\[
\la_i \| \Om \|_\ell \leq C_2 \| \Om \|_{\ell +1} 
+ \| h \|_\ell \ .
\]
Then we have
\[
\la_i^2 \| \Om \|^2_\ell \leq C_3 \| \Om \|^2_{\ell +1}
+ C_4 \| h \|^2_\ell \ ,
\]
that is
\begin{equation}
\frac{\nu}{2C_3}\la_i^2 \| \Om \|^2_\ell \leq \frac{1}{2}\nu 
\| \Om \|^2_{\ell +1} + \frac{\nu C_4}{2C_3} \| h \|^2_\ell \ .
\label{SDS3}
\end{equation}
Adding (\ref{SDS2}) and (\ref{SDS3}), we have
\[
\left ( \la_r + \frac{\nu}{2C_3}\la_i^2 -C_1 \right ) \| \Om \|^2_\ell \leq  
\frac{1}{2} \left (1+ \frac{\nu C_4}{C_3} \right ) \| h \|^2_\ell \ .
\]
Let $a = C_1$, $b=\frac{\nu}{2C_3}$, and $C_5 =\frac{1}{2} \left ( 
1+ \frac{\nu C_4}{C_3} \right )$,
\[
( \la_r -a +b \la_i^2) \| \Om \|^2_\ell \leq C_5 \| h \|^2_\ell \ .
\]
Thus if $\la_r -a +b \la_i^2 > 0$,
\begin{equation}
\| (L-\la )^{-1} \| \leq \frac{C_5^{1/2}}{( \la_r -a +b \la_i^2)^{1/2}} \ .
\label{SDS4}
\end{equation}
$(L-\la )^{-1}$ is holomorphic in $\la$ and can be continued analytically 
as long as $\| (L-\la )^{-1} \|$ is bounded. As shown above, $(L-\la )^{-1}$
exists at large $\la > 0$. Then from (\ref{SDS4}), $(L-\la )^{-1}$
exists everywhere outside the parabolic region $\la_r  \leq a-b \la_i^2$.
Thus the eigenvalues of $L$ lie inside this region. To prove the rest of the 
theorem, we need a known lemma [Corollary 31, pp.1042, \cite{DS63}].
\begin{lemma}
Let $T$ be a densely defined unbounded linear operator in a Hilbert space 
$H$, with the property that for some $\la_0$ in the resolvent set of $T$,
$(T-\la_0 )^{-1}$ is a Hilbert-Schmidt operator. Let $\ga_1, \cdots, \ga_5$
be non-overlapping differentiable arcs having limiting directions at 
infinity, and such that no adjacent pair of arcs form an angle as great as 
$\pi /2$ at infinity. Suppose that the resolvent $(T-\la )^{-1}$ satisfies 
an inequality $\| (T-\la )^{-1} \| = O(|\la |^{-N})$ as $\la \ra \infty$ 
along each arc $\ga_j$ for some natural number $N$. Then the eigenfunctions 
(and higher order eigenfunctions) form a complete basis in $H$.
\label{Dun}
\end{lemma}
\nid
From (\ref{SDS4}), we see that along any ray $\la_r = c \la_i$ ($c$ finite),
\[
\| (L-\la )^{-1} \| = O(|\la |^{-1}) \quad \mbox{as} \quad \la \ra \infty \ .
\]
Next we show that $(L-\la )^{-1}$ is a Hilbert-Schmidt operator for large 
enough $\la > 0$. Notice that 
\[
(L-\la )^{-1} = \left [ I +(\nu \Dl -\la )^{-1}A \right ]^{-1}
(\nu \Dl -\la )^{-1}\ .
\]
Since
\begin{eqnarray*}
& & \| (\nu \Dl -\la )^{-1} \Om \|^2_\ell = \sum_{k} \frac{|k|^{2\ell }}
{(\nu |k|^2 +\la)^2}|\Om_k|^2 \\
& & = \frac{1}{\nu^2} \sum_{k} \frac{1}{|k|^2 +\la / \nu}\frac{ |k|^2}
{|k|^2 +\la / \nu} |k|^{2(\ell -1)}|\Om_k|^2 \\
& & \leq \frac{1}{\nu \la }\sum_{k} |k|^{2(\ell -1)}|\Om_k|^2 \ ,
\end{eqnarray*}
we have 
\[
\| (\nu \Dl -\la )^{-1} \Om \|_\ell \leq \frac{1}{\sqrt{\nu \la}} 
\| \Om \|_{\ell -1} \ .
\]
It is easy to see that
\[
\| A\Om \|_{\ell -1} \leq C \| \Om \|_\ell \ .
\]
Thus
\[
\| (\nu \Dl -\la )^{-1}A \Om \|_\ell \leq \frac{C}{\sqrt{\nu \la}} 
\| \Om \|_\ell \ .
\]
Choose $\la_0 = \frac{2C^2}{\nu}$, then
\[
\| (\nu \Dl -\la_0 )^{-1}A \| \leq \frac{1}{2} \ .
\]
We have $\left [ I +(\nu \Dl -\la_0 )^{-1}A \right ]^{-1}$ being a bounded 
operator. Recall that a bounded linear operator $T$ is a Hilbert-Schmidt
operator if
\[
\sum \| Te_k\|^2
\]
is finite where $\{ e_k\}$ is a complete orthonormal basis \cite{DS63}.
By using the Fourier basis for $H^\ell$, one sees that 
$(\nu \Dl -\la_0 )^{-1}$ is a Hilbert-Schmidt operator. A product of a 
bounded linear operator with a Hilbert-Schmidt operator is another 
Hilbert-Schmidt operator. Thus $(L -\la_0 )^{-1}$ is a Hilbert-Schmidt 
operator. By Lemma \ref{Dun}, the corresponding eigenfunctions (and 
higher order eigenfunctions) of $L$ form a complete basis in $H^\ell$ 
($\ell = 0,1,2, \cdots$).
\end{proof}
\nid
We will study the case as described by the following set-up.
\begin{itemize}
\item {\bf Set-Up}. Let there be $m$ unstable eigenvalues $\la_1^u 
\cdots \la_m^u$, $n$ center eigenvalues $\la_1^c \cdots \la_n^c$, 
and of course the rest infinitely many stable eigenvalues 
$\la_1^s \ \la_2^s \cdots $. Let the real parts of $\la_1^u$ and $\la_1^s$ 
have the smallest absolute values among the unstable and stable 
eigenvalues, respectively. 
\end{itemize}
\begin{remark}
For specific examples, 
often there is no center eigenvalue, here we put them in for generality.
We need 
$H^\ell$ ($\ell =0,1,2,\cdots$) to be a Banach algebra \cite{Ada75}. 
Thus for 2D and 3D, $\ell \geq 2$. Our argument requires that we work 
with $H^\ell$ where $\ell \geq 3$ for both 2D and 3D.
\end{remark}
We have the splitting
\[
H^\ell = H^\ell_u + H^\ell_c + H^\ell_s
\]
where $H^\ell_z$ ($z=u,c,s$) are spanned by the corresponding eigenfunctions.
Projections onto  $H^\ell_z$ are denoted by $P^z$ ($z=u,c,s$). Denote by 
\[
L^z = P^z L \ , \quad N^z = P^z N \ , \quad \Om^z = P^z \Om \ , \quad 
(z=u,c,s)\ .
\]
Then (\ref{NSr}) can be rewritten as
\begin{equation}
\pa_t \Om^z = L^z \Om^z +N^z(\Om)\ , \quad (z=u,c,s)\ . 
\label{NSrr}
\end{equation}
Since in this case the spectral mapping theorem is trivially true, one has 
the trichotomy:
\begin{eqnarray}
& & \| e^{tL^u} \| \leq C e^{\al t} \ , \quad t \leq 0\ ; \label{triu} \\
& & \| e^{tL^c} \| \leq C |t|^{n_1} \ , \quad t \in \mathbb{R}\ ; 
\label{tric} \\
& & \| e^{tL^s} \| \leq C e^{-\be t} \ , \quad t \geq 0\ ; \label{tris}
\end{eqnarray}
where $\al = \ \mbox{Re}\{ \la_1^u \} -\e$, $\be = \ -\mbox{Re}\{ \la_1^s 
\} -\e$, $ n_1 \leq n$, and $0 < \e \ll  -\mbox{Re}\{ \la_1^s \}$. Our 
argument depends upon the following two crucial facts:
\begin{enumerate}
\item One has a stronger inequality than (\ref{tris}) [\cite{Hen81}, 
Theorem 1.5.4] due to the fact that $L^s$ is sectorial,
\begin{equation}
\| e^{tL^s}\Om^s \|_{\ell +1} \leq C(\be t)^{-1/2} e^{-\be t} 
\| \Om^s\|_\ell \ , \quad t \geq 0 \ . 
\label{tris1}
\end{equation}
\item Due to finite dimensionality, one has 
\begin{equation}
\| \Om^z \|_{\ell +1} \leq C \| \Om^z\|_\ell \ , \quad (z=u,c)\ . 
\label{bdd}
\end{equation}
\end{enumerate}
\begin{remark}
For (\ref{bdd}), take $z=u$ for example, since $\Om_1^u \cdots \Om_m^u$
belong to $H^\ell$ for any $\ell$, we can take  
\begin{equation}
C=\max_{1\leq j \leq m} \left \{ \| \Om_j^u \|_{\ell +1}\bigg / 
\| \Om_j^u \|_\ell \right \} \ .
\label{bdd1}
\end{equation}
\end{remark}
\begin{theorem}
In a neighborhood of the fixed point $\Om_*$ in the Sobolev space 
$H^\ell (\mathbb{T}^d)$ ($\ell \geq 3, d=2,3$), there exist a 
$m$-dimensional $C^\infty$ unstable manifold, a $n$-dimensional $C^\infty$
center manifold, a ($m+n$)-codimensional $C^\infty$ stable manifold, 
a ($m+n$)-dimensional $C^\infty$ center-unstable manifold, and a 
$m$-codimensional $C^\infty$ center-stable manifold.
\label{MTHM}
\end{theorem}
\begin{proof}
We will take the center-stable manifold as the example.
The arguments for the others are similar.
Apply the method of variation of parameters to (\ref{NSrr}),
one gets the integral equations
\begin{equation}
\Om^z(t) = e^{(t-t_0)L^z}\Om^z(t_0)+\int_{t_0}^t e^{(t-\tau)L^z}
N^z(\Om)d\tau \ , \quad (z=u,c,s)\ .
\label{inte}
\end{equation}
To prove the existence of a center-stable manifold, we start with the Banach 
space 
\[
B_{\ga,\ell}=\left \{ \Om(t)\ \bigg |\ \| \Om \|_{\ga,\ell} = 
\sup_{t \geq 0} e^{-\ga t} \sum_{z=u,c,s}\| \Om^z(t) \|_\ell < +\infty\ ,
\quad (\frac{1}{4} \al < \ga < \al )\right \}\ .
\]
For such $\Om(t)$ in $B_{\ga,\ell}$
\[
\lim_{t_0 \ra +\infty}e^{(t-t_0)L^u}\Om^u(t_0) \ra 0 \ .
\]
In the equation (\ref{inte}), for $z=u$, let $t_0 \ra +\infty$; 
for $z=c,s$, let $t_0=0$; then one gets
\begin{eqnarray}
\Om^u(t) &=& \int_{+\infty}^t e^{(t-\tau)L^u}
N^u(\Om)d\tau \ , \label{inte1} \\
\Om^z(t) &=& e^{tL^z}\Om^z(0)+\int_{0}^t e^{(t-\tau)L^z}
N^z(\Om)d\tau \ , \quad (z=c,s)\ . \label{inte2}
\end{eqnarray}
The right hand side of (\ref{inte1})-(\ref{inte2}) defines a map $\Ga =
\Ga^u + \Ga^c + \Ga^s$ on $B_{\ga,\ell}$. First we show that $\Ga$ contracts,
then $\Ga$ acting is an easier argument. From (\ref{inte1}), one has 
\[
e^{-\ga t} \| \Ga^u \Om_1 - \Ga^u \Om_2 \|_\ell \leq C \int^{+\infty}_t
e^{-\ga t}e^{\al (t-\tau)} \| N^u(\Om_1) - N^u(\Om_2) \|_\ell d \tau \ .
\]
Using (\ref{bdd}), one gets
\[
e^{-\ga t} \| \Ga^u \Om_1 - \Ga^u \Om_2 \|_\ell \leq C \int^{+\infty}_t
e^{-\ga t}e^{\al (t-\tau)} \| N^u(\Om_1) - N^u(\Om_2) \|_{\ell -1} d \tau \ ,
\]
in this article, all constants are denoted by the same $C$. Since 
$H^{\ell -1}$ is a Banach algebra
\[
\| N^u(\Om_1) - N^u(\Om_2) \|_{\ell -1} \leq C(\| \Om_1 \|_\ell 
+\| \Om_2 \|_\ell )\| \Om_1 - \Om_2 \|_\ell \ .
\]
Finally one gets
\[
e^{-\ga t} \| \Ga^u \Om_1 - \Ga^u \Om_2 \|_\ell \leq C \int^{+\infty}_t
e^{(\al -\ga )(t-\tau)}(\| \Om_1 \|_\ell +\| \Om_2 \|_\ell )e^{-\ga \tau }
\| \Om_1 - \Om_2 \|_\ell d \tau \ .
\]
Now we need to replace $N(\Om)$ by its cut-off $\chi (\| \Om \|_\ell /\dl)
N(\Om)$ for some $\dl > 0$, where
\begin{eqnarray}
& & \chi (r) = 1\ , \quad r \in [0,1]\ ; \quad 
\chi (r) = 0\ , \quad r \in [3, +\infty )\ ; \label{cutoff} \\
& & \chi'(r) < 1 \ , \quad r \in [0, +\infty )\ ;
\quad \chi \in C_0^\infty (\mathbb{R}, \mathbb{R})\ . \nonumber
\end{eqnarray}
If both $\chi (\| \Om_1 \|_\ell /\dl)$ and $\chi (\| \Om_2 \|_\ell /\dl)$
are zero, then of course
\[
\chi (\| \Om_1 \|_\ell /\dl)N(\Om_1) - \chi (\| \Om_2 \|_\ell /\dl)N(\Om_2)
=0\ .
\]
If one of them is nonzero, without loss of generality, say $\chi 
(\| \Om_1 \|_\ell /\dl) \neq 0$, then
\begin{eqnarray*}
& & \| \chi (\| \Om_1 \|_\ell /\dl)N(\Om_1) - \chi (\| \Om_2 \|_\ell /\dl)
N(\Om_2) \|_{\ell -1} \\
&\leq& \| [\chi (\| \Om_1 \|_\ell /\dl) - \chi (\| \Om_2 \|_\ell /\dl)]
N(\Om_1) \|_{\ell -1} + \| \chi (\| \Om_2 \|_\ell /\dl) [N(\Om_1) -
N(\Om_2)]\|_{\ell -1} \\
&\leq& | \chi (\| \Om_1 \|_\ell /\dl) - \chi (\| \Om_2 \|_\ell /\dl)| 
\| N(\Om_1) \|_{\ell -1} + |\chi (\| \Om_2 \|_\ell /\dl)| \| N(\Om_1) -
N(\Om_2)\|_{\ell -1} \\
&\leq& \frac{1}{\dl} \| \Om_1 - \Om_2 \|_\ell \| N(\Om_1)\|_{\ell -1} 
+|\chi (\| \Om_2 \|_\ell /\dl)| \| N(\Om_1) - N(\Om_2)\|_{\ell -1} \ .
\end{eqnarray*}
If $\chi (\| \Om_2 \|_\ell /\dl) = 0$, then the second part disappears; 
otherwise, $\| \Om_2 \|_\ell < 3\dl$. In any case, one gets
\[
\| \chi (\| \Om_1 \|_\ell /\dl)N(\Om_1) - \chi (\| \Om_2 \|_\ell /\dl)
N(\Om_2) \|_{\ell -1} \leq C \dl \| \Om_1 - \Om_2 \|_\ell \ .
\]
Thus, one has
\begin{eqnarray*}
e^{-\ga t} \| \Ga^u \Om_1 - \Ga^u \Om_2 \|_\ell 
&\leq& C\dl \| \Om_1 - \Om_2 \|_{\ga,\ell} \int^{+\infty}_t
e^{(\al -\ga )(t-\tau)}d \tau \\
&=& C\frac{\dl}{\al -\ga }\| \Om_1 - \Om_2 \|_{\ga,\ell} \ . 
\end{eqnarray*}
Thus
\begin{equation}
\sup_{t \geq 0} e^{-\ga t} \| \Ga^u \Om_1 - \Ga^u \Om_2 \|_\ell  
\leq C\frac{\dl}{\al -\ga }\| \Om_1 - \Om_2 \|_{\ga,\ell} \ .
\label{est1}
\end{equation}
For $z=c$ in (\ref{inte2}), one has
\[
e^{-\ga t} \| \Ga^c \Om_1 - \Ga^c \Om_2 \|_\ell \leq 
C \int^t_0 e^{-\ga t} e^{\e (t-\tau)}\| N^c(\Om_1) -
N^c(\Om_2)\|_\ell d \tau \ ,
\]
here notice that $\Om_1^c$ and $\Om_2^c$ have the same initial 
condition $\Om_1^c(0)=\Om_2^c(0)$. The same argument as above leads to 
\begin{eqnarray*}
e^{-\ga t} \| \Ga^c \Om_1 - \Ga^c \Om_2 \|_\ell 
&\leq& C\dl \| \Om_1 - \Om_2 \|_{\ga,\ell} \int_0^t
e^{-(\ga -\e )(t-\tau)}d \tau \\
&\leq& C\frac{\dl}{\ga -\e}\| \Om_1 - \Om_2 \|_{\ga,\ell} \ . 
\end{eqnarray*}
Thus
\begin{equation}
\sup_{t \geq 0} e^{-\ga t} \| \Ga^c \Om_1 - \Ga^c \Om_2 \|_\ell  
\leq C\frac{\dl}{\ga -\e }\| \Om_1 - \Om_2 \|_{\ga,\ell} \ .
\label{est2}
\end{equation}
For $z=s$ in (\ref{inte2}), one has
\[
e^{-\ga t} \| \Ga^s \Om_1 - \Ga^s \Om_2 \|_\ell \leq 
C \int^t_0 e^{-\ga t} [\be (t-\tau )]^{-1/2}e^{-\be (t-\tau)}\| N^s(\Om_1) -
N^s(\Om_2)\|_{\ell -1} d \tau \ .
\]
The same argument as before leads to 
\begin{eqnarray*}
e^{-\ga t} \| \Ga^s \Om_1 - \Ga^s \Om_2 \|_\ell 
&\leq& C\dl \| \Om_1 - \Om_2 \|_{\ga,\ell} \int_0^t
[\be (t-\tau )]^{-1/2}e^{-(\be +\ga )(t-\tau)}d \tau \\
&\leq& C\frac{\dl}{\sqrt{\be}}\| \Om_1 - \Om_2 \|_{\ga,\ell} \ . 
\end{eqnarray*}
Thus
\begin{equation}
\sup_{t \geq 0} e^{-\ga t} \| \Ga^s \Om_1 - \Ga^s \Om_2 \|_\ell  
\leq C\frac{\dl}{\sqrt{\be} }\| \Om_1 - \Om_2 \|_{\ga,\ell} \ .
\label{est3}
\end{equation}
By choosing 
\begin{equation}
\dl = \min \left \{ \frac{1}{6} C(\al -\ga ), \frac{1}{6} C(\ga -\e ),
\frac{1}{6} C\sqrt{\be} \right \}\ ,
\label{chos}
\end{equation}
(\ref{est1})-(\ref{est3}) imply that
\[
\| \Ga \Om_1 - \Ga \Om_2 \|_{\ga,\ell} \leq \frac{1}{2}  
\| \Om_1 - \Om_2 \|_{\ga,\ell} \ .
\]
Thus $\Ga$ contracts. The claim that $\Ga$ acts, i.e. $\Ga : 
B_{\ga,\ell} \ra B_{\ga,\ell}$; follows similarly, and is an easier argument. 
By the contraction mapping theorem, $\Ga$ has a fixed point $\Om$ in 
$B_{\ga,\ell}$, which satisfies (\ref{inte1})-(\ref{inte2}), of course, with 
$N$ replaced by its cut-off. Equations (\ref{inte1})-(\ref{inte2}) define 
a map which maps $\Om^z(0)$ ($z=c,s$) into 
\[
\Om^u(0)=\int_{+\infty}^0 e^{-\tau L^u} N^u(\Om) d \tau \ .
\]
This map defines the center-stable manifold. Due to the cut-off, only 
in the $\dl$ neighborhood of $\Om_*$, this center-stable 
manifold corresponds to the NS flow. Smoothness can be proved through 
the standard argument \cite{Li04}. This proves the theorem. 
\end{proof}

To study the zero-viscosity limits of invariant manifolds, we need to make 
an assumption.
\begin{itemize}
\item {\bf Assumption 1}. As $\nu \ra 0^+$, $\mbox{Re}\{\la_1^u\}$ is of 
order $O(1)$, and $\mbox{Re}\{\la_1^s\}$ is $O(\nu)$, and the constant 
$C$ in (\ref{triu})-(\ref{bdd}) is $O(1)$.
\end{itemize}
\begin{theorem}
Under Assumption 1, as $\nu \ra 0^+$, in an order $O(\sqrt{\nu})$ 
neighborhood of the fixed point $\Om_*$ in the Sobolev space 
$H^\ell (\mathbb{T}^d)$ ($\ell \geq 3, d=2,3$), there exist an unstable 
manifold and a center-stable manifold; and in its order $O(\nu )$ 
neighborhood, there exist a stable manifold, a center manifold, and 
a center-unstable manifold.
\end{theorem}
\begin{proof}
Under the Assumption 1, the $\dl$ in (\ref{chos}) is of order $O(\sqrt{\nu})$
since $\be$ is $O(\nu)$. Thus we have the claim for the center-stable manifold.
By keeping track of similar estimates as in the proof of Theorem \ref{MTHM},
one can quickly see the rest of the claims.
\end{proof}
Due to Assumption 1, zero-viscosity limits of invariant manifolds become 
a challenging problem. The challenge is further amplified by the elusiveness
of the invariant manifolds for Euler equations as discussed in the 
Introduction. As a test, we studied a simpler problem in \cite{Li05}, 
which reveals some unique features of the zero-viscosity limit of the 
spectra. For instance, the zero-viscosity limit of a discrete spectrum 
can be a continuous spectrum which is not the spectrum of zero-viscosity.

\section{Examples}

Using Fourier series for the 2D NS (\ref{2DNS}),
\[
\Om = \sum_{k \in \ZZo} \om_k e^{ik\cdot x}\ , \quad 
f = \sum_{k \in \ZZo} f_k e^{ik\cdot x}\ , 
\]
where $\om_{-k} = \overline{\om_k}$ and $f_{-k} = \overline{f_k}$, one gets
the kinetic form of 2D NS
\[
\dot{\om}_k = \sum_{k=q+r} A(q,r) \ \om_q \om_r +\nu [-|k|^2 \om_k +f_k ] \ ,
\]
where 
\[
A(q,r) = \frac{1}{2}\left [|r|^{-2}-|q|^{-2}\right ]\left | \begin{array}{lr} 
q_1 & r_1 \\ q_2 & r_2 \\ \end{array} \right | \ , 
\]
where $|q|^2 =q_1^2 +q_2^2$ for $q=(q_1,q_2)^T$, similarly for $r$.
Linearize the 2D NS at the steady state given by a single mode,
\begin{equation}
\Om_* = \Ga e^{ip\cdot x} + \bar{\Ga} e^{-ip\cdot x}
\label{sm}
\end{equation}
where $f = |p|^2 \Om_*$, one gets
\begin{eqnarray}
\dot{\omega}_{\hat{k} + np} &=& A(p, \hat{k} + (n-1) p) 
     \ \Gamma \ \omega_{\hat{k} + (n-1) p} \nonumber \\  
&+& \ A(-p, \hat{k} + (n+1)p)\ 
     \bar{\Gamma} \ \omega_{\hat{k} +(n+1)p} -\nu | \hk +np|^2 
\om_{\hk +np}\ , \label{decu}
\end{eqnarray}
where its right hand side defines the linear NS operator $L$. Let
$\om_{\hk +np} = e^{\la t}\tom_{\hk +np}$, $\Ga = |\Ga| e^{i\ga}$, 
$a = \frac{1}{2} |\Ga| \left | \begin{array}
{lr} p_1 & \hat{k}_1 \\ p_2 & \hat{k}_2 \\ \end{array} \right |$, 
$\rho_n = |p|^{-2} - | \hat{k}+np|^{-2}$, and $z_n =\rho_n 
e^{-in\ga}\tom_{\hk +np}$, then
\begin{equation}
a_n z_n +z_{n-1} - z_{n+1} = 0 \ , 
\label{PP}
\end{equation}
where $a_n = (a \rho_n)^{-1}(\la + \nu | \hk +np|^2)$. Let 
$w_n = z_n / z_{n-1}$, one gets \cite{MS61} 
\begin{equation}
a_n + \frac{1}{w_n} = w_{n+1}\ ,
\label{IR1}
\end{equation}
which leads to the continued fraction,
\begin{equation}
w_n^{(1)}=a_{n-1} +\frac{1}{a_{n-2} + \frac{1}{a_{n-3}+\cdots}}\ \ .
\label{CF1}
\end{equation}
Rewriting (\ref{IR1}) as,
\begin{equation}
w_n = \frac{1}{-a_n +w_{n+1}}\ , 
\label{IR2}
\end{equation}
which leads to the other continued fraction 
\begin{equation}
w_n^{(2)}=-\frac{1}{a_{n} + \frac{1}{a_{n+1}+\cdots}}\ \ .
\label{CF2}
\end{equation}   
Along the spirit of the finite difference of a second order ordinary 
differential equation, the difference equation (\ref{PP}) should have two 
linearly independent solutions. Rigorous theory has been well developed. 
This is the so-called Poincar\'e-Perron system. For details, see \cite{Poi85}
\cite{Per10a} \cite{Per10b} \cite{Per11} \cite{MT33} \cite{Gau67} \cite{JT80}
\cite{LW92} \cite{Li00} \cite{Liu95}. When $\nu >0$, $a_n \ra \ta n^2$ 
as $|n| \ra +\infty$. Then there are a growing solution $z_n /z_{n-1} \ra 
\ta n^2$ and a decaying solution $z_n /z_{n-1} \ra \ta^{-1} n^{-2}$, as 
$n \ra +\infty$; and vice versa as $n \ra -\infty$. The intuition on this 
is clear from (\ref{IR1}) and (\ref{IR2}). The eigenvalues are then determined 
by matching the two decaying solutions given by the two continued fractions 
(\ref{CF1}) and (\ref{CF2}). When $\nu = 0$, $a_n \ra \ta = a |p|^2 \la$ as 
$|n| \ra +\infty$. When $\mbox{Re}\{ \ta \} \neq 0$, or 
$\mbox{Re}\{ \ta \}=0 \ (|\ta| >2)$, there are a growing solution 
$z_n /z_{n-1} \ra w_+$, ($|w_+| > 1$), and a decaying solution 
$z_n /z_{n-1} \ra w_-$, ($|w_-| < 1$), as $n \ra +\infty$; and vice versa 
as $n \ra -\infty$. Here $w_\pm$ solve the characteristic equation 
\cite{Li00}
\[
w^2 -\ta w - 1 = 0 \ .
\]
Again the eigenvalues are determined by matching the two decaying solutions 
given by the two continued fractions (\ref{CF1}) and (\ref{CF2}). Finally,
$\mbox{Re}\{ \ta \}=0 \ (|\ta| \leq 2)$ corresponds to a continuous spectrum
\cite{Li00}. The proof of \cite{Li00} can be generalized to $H^\ell$ setting 
for any natural number $\ell$. For more general result, see \cite{SL03}. In 
summary, when $\nu >0$ or $\nu = 0$ [$\mbox{Re}\{ \ta \} \neq 0$, or 
$\mbox{Re}\{ \ta \}=0 \ (|\ta| >2)$], the eigenvalues are determined by
\begin{equation}
a_{0} +\frac{1}{a_{-1} + \frac{1}{a_{-2}+\cdots}} =
-\frac{1}{a_{1} + \frac{1}{a_{2}+\cdots}} \ ,
\label{mch}
\end{equation}
which is obtained by matching $w^{(1)}_n$ (\ref{CF1}) and $w^{(2)}_n$ 
(\ref{CF2}) at $n=1$, i.e. $w^{(1)}_1=w^{(2)}_1$. The corresponding 
eigenfunctions belong to $H^\ell$ for any natural number $\ell$. 

\subsection{Example 1}

Our first prime example is the shear $\Om_* = \cos x_2$ defined on the 
rectangular periodic domain $[0, 2\pi /\al ] \times [0, 2\pi ]$ where 
$1/2 < \al < 1$. As shown below, this example has only one unstable 
eigenvalue. Complete information on this eigenvalue can be obtained. 
Thereby, existence of an unstable manifold and a stable manifold 
can be confirmed. This example is motivated by examples 2 and 3 studied 
later.

In this rectangular domain case, $\Om$ has the Fourier expansion
\[
\Om =\sum_{k \in \ZZo} \om_k e^{i[\al k_1 x_1 +k_2 x_2]}\ .
\]
We have equation (\ref{PP}) with 
\[
a_n = - \frac{4}{\al \hk_1}\frac{(\al \hk_1)^2+(\hk_2+n)^2}
{(\al \hk_1)^2+(\hk_2+n)^2-1} \left \{ \la +\nu \left [ 
(\al \hk_1)^2+(\hk_2+n)^2 \right ] \right \} \ .
\]
\begin{theorem}
The spectra of the 2D linear Euler operator $L$ have the following 
properties.
\begin{enumerate}
\item $(\al \hk_1)^2+(\hk_2+n)^2 > 1$, $\forall n \in \Z$. 
When $\nu \ra 0$, there is no eigenvalue of non-negative real part. 
When $\nu = 0$, the entire spectrum is the continuous spectrum
\[
\left [ -i\frac{\al |\hk_1|}{2}, \ i\frac{\al |\hk_1|}{2} \right ]\ .
\]
\item $\hk_1 = 0$, $\hk_2 = 1$. The spectrum consists of the eigenvalues 
\[
\la = - \nu n^2 \ , \quad n \in \Zo \ .
\]
The eigenfunctions are the Fourier modes
\[
\tom_{np} e^{inx_2} + \ \mbox{c.c.}\ \ , \quad \forall \tom_{np} \in 
\C\ , \quad n \in \Zo \ .
\]
As $\nu \ra 0^+$, the eigenvalues are dense on the negative half of the real 
axis $(-\infty, 0]$. Setting $\nu =0$, the only eigenvalue is $\la = 0$ of 
infinite multiplicity with the same eigenfunctions as above.
\item $\hk_1 = -1$, $\hk_2 = 0$. (a). $\nu >0$. For any $\al \in (0.5, 0.95)$,
there is a unique $\nu_*(\al)$,
\begin{equation}
\frac{\sqrt{32-3\al^6-17\al^4-16\al^2}}{4(\al^2+1)(\al^2+4)} < 
\nu_*(\al) < \frac{1}{2(\al^2+1)} \sqrt{\frac{1-\al^2}{2}}\ ,
\label{nuda}
\end{equation}
where the term under the square root on the left is positive for 
$\al \in (0.5, 0.95)$, and the left term is always less than the right term.
When $\nu > \nu_*(\al)$, there is no eigenvalue of non-negative real part. 
When $\nu = \nu_*(\al)$, $\la =0$ is an eigenvalue, and all the rest 
eigenvalues have negative real parts. When $\nu < \nu_*(\al)$, there is 
a unique positive eigenvalue $\la (\nu )>0$, and all the rest 
eigenvalues have negative real parts. $\nu^{-1} \la (\nu )$ is a strictly 
monotonically decreasing function of $\nu$. When $\al \in (0.5, 0.8469)$,
we have the estimate
\begin{eqnarray*}
& & \sqrt{\frac{\al^2(1-\al^2)}{8(\al^2+1)}-\frac{\al^4 (\al^2+3)}{16
(\al^2+1)(\al^2+4)}} - \nu (\al^2+1) < \la (\nu ) \\
& & < \sqrt{\frac{\al^2(1-\al^2)}{8(\al^2+1)}}- \nu \al^2 \ ,
\end{eqnarray*}
where the term under the square root on the left is positive for 
$\al \in (0.5, 0.8469)$.
\[
\sqrt{\frac{\al^2(1-\al^2)}{8(\al^2+1)}-\frac{\al^4 (\al^2+3)}{16
(\al^2+1)(\al^2+4)}} \leq \lim_{\nu \ra 0^+} \la (\nu )  \leq 
\sqrt{\frac{\al^2(1-\al^2)}{8(\al^2+1)}} \ .
\]
In particular, as $\nu \ra 0^+$, $\la (\nu ) =O(1)$.

(b). $\nu =0$. When $\al \in (0.5, 0.8469)$, we have only two eigenvalues
$\la_0$ and $-\la_0$, where $\la_0$ is positive,
\[
\sqrt{\frac{\al^2(1-\al^2)}{8(\al^2+1)}-\frac{\al^4 (\al^2+3)}{16
(\al^2+1)(\al^2+4)}} < \la_0 <
\sqrt{\frac{\al^2(1-\al^2)}{8(\al^2+1)}} \ .
\]
The rest of the spectrum is a continuous spectrum $[-i\al /2, \ i\al /2]$.

(c). For any fixed $\al \in (0.5, 0.8469)$,
\begin{equation}
\lim_{\nu \ra 0^+} \la (\nu ) = \la_0 \ .
\label{pet1}
\end{equation}
\item Finally, when $\nu = 0$, the union of all the above pieces of 
continuous spectra is the imaginary axis $i\R$.
\end{enumerate}
\label{PET}
\end{theorem}
\nid
From Theorems \ref{MTHM} and \ref{PET}, we have the corollary.
\begin{corollary}
For any $\al \in (0.5, 0.95)$, and $\nu \in (0, \nu_*(\al ))$ where 
$\nu_*(\al ) > 0$ satisfies (\ref{nuda}), in a neighborhood of $\Om_*$
in the Sobolev space $H^\ell (\mathbb{T}^2)$ ($\ell \geq 3$),
there are an $1$-dimensional $C^\infty$ unstable manifold and an 
$1$-codimensional $C^\infty$ stable manifold.
\end{corollary}
\begin{remark}
In the Theorem \ref{PET}, (\ref{pet1}) verifies part of Assumption 1. 
That is, $\mbox{Re}\{\la_1^u\}$ is of order $O(1)$ as $\nu \ra 0^+$. 
Case 2 of Theorem \ref{PET} indicates that $\mbox{Re}\{\la_1^s\}$ is
at least $O(\nu)$ as $\nu \ra 0^+$. Case 3 indicates that the constant 
$C$ in (\ref{triu}) and (\ref{bdd}) should be $O(1)$ as $\nu \ra 0^+$.
\end{remark}
\begin{proof}
We will give the proof case by case.
\begin{enumerate}
\item The case (a). $\nu > 0$. If $\mbox{Re}\{ \la \} \geq 0$, then all 
the $\mbox{Re}\{ a_n \}$'s are nonzero and have the same sign. The real 
parts of the right and left hand sides of (\ref{mch}) are nonzero but of 
different signs. Thus there is no eigenvalue of non-negative real part. 
The case (b). $\nu = 0$. If $\mbox{Re}\{ \la \} \neq 0$, then all 
the $\mbox{Re}\{ a_n \}$'s are nonzero and have the same sign. Again 
(\ref{mch}) can not be satisfied. If $\mbox{Re}\{ \la \} = 0$, let 
$\ta = \lim_{n \ra \infty} a_n = -\frac{4}{\al \hk_1} \la$, then
\[
\ta z_n = \frac {(\al \hk_1)^2 + (\hk_2 +n)^2 - 1}
{(\al \hk_1)^2 + (\hk_2 +n)^2}(z_{n+1}-z_{n-1})\ ,
\]
where $(\al \hk_1)^2 + (\hk_2 +n)^2 > 1$, $\forall n \in \Z$. By using 
$\ell_2$ norm of $\{ z_n \}_{n \in \Z}$, one sees that possible 
eigenvalues have to satisfy $|\ta| \leq 2$. As mentioned before, 
$\mbox{Re}\{ \la \} = 0$ and $|\ta | \leq 2$ correspond to a continuous spectrum \cite{Li00}, which is the interval
\begin{equation}
\left [ -i\frac{\al |\hk_1|}{2}, \ i\frac{\al |\hk_1|}{2} \right ]\ .
\label{pie}
\end{equation}
\item In this case, one has
\[
[\la +\nu (n+1)^2]\tom_{(n+1)p} = 0 \ .
\]
The claims follow immediately.
\item In this case,
\[
a_n = \frac{4}{\al} \frac{\al^2+n^2}{\al^2+n^2-1}[\la +\nu (\al^2+n^2)]\ .
\]
Thus $a_{-n}=a_n$. Equation (\ref{mch}) is reduced to 
\begin{equation}
-a_{0}/2 = \frac{1}{a_{1} + \frac{1}{a_{2}+\cdots}} \ .
\label{prd1}
\end{equation}
The first a few $a_n$'s are 
\begin{eqnarray*}
a_0 &=& \frac{4}{\al} \frac{\al^2}{\al^2-1}[\la +\nu \al^2]\ , \\
a_1 &=& \frac{4}{\al} \frac{\al^2+1}{\al^2}[\la +\nu (\al^2+1)]\ , \\
a_2 &=& \frac{4}{\al} \frac{\al^2+4}{\al^2+3}[\la +\nu (\al^2+4)]\ . 
\end{eqnarray*}
Let 
\[
f(\la ) = \frac{1}{a_{1} + \frac{1}{a_{2}+\cdots}} \ , \quad 
g(\la ) = -a_{0}/2 \ .
\]
(a). $\nu >0$. First we show that if $\la$ ($\mbox{Re}\{ \la \} \geq 0$) 
is a solution of (\ref{prd1}), then $\la$ must be real. Assume  
$\mbox{Im}\{ \la \} > 0$, then
\[
\arg (-a_0) > \arg (a_1) > \arg (a_2) > \cdots \geq 0 \ .
\]
Thus 
\[
|\arg (f(\la ))| \leq \arg (a_1) \ .
\]
But 
\[
\arg (g(\la )) = \arg (-a_0) > |\arg (f(\la ))| \ ,
\]
a contradiction. Similarly for the case $\mbox{Im}\{ \la \} < 0$, thus 
$\la$ ($\mbox{Re}\{ \la \} \geq 0$) is real. Next we show that there is 
a $\nu_0 >0$ such that for every $\nu \leq \nu_0$, there is a unique 
eigenvalue $\la >0$. Since $a_n >0$, $\forall n \geq 1$, we have
\begin{equation}
\frac{1}{a_{1} + \frac{1}{a_{2}}} < f(\la ) < \frac{1}{a_{1}}\ .
\label{prd2}
\end{equation}
Thus when $\la$ is large enough
\begin{equation}
f(\la ) < g(\la ) \ .
\label{prd3}
\end{equation}
From (\ref{prd2}), 
\[
f(0) > \frac{1}{\frac{4\nu (\al^2+1)^2}{\al^3} + \frac{\al (\al^2+3)}
{4\nu (\al^2+4)^2}}\ .
\]
We know that 
\[
g(0) =  \frac{ 2 \nu \al^3}{1-\al^2} \ .
\]
We want to choose $\nu$ such that 
\[
\frac{1}{\frac{4\nu (\al^2+1)^2}{\al^3} + \frac{\al (\al^2+3)}
{4\nu (\al^2+4)^2}} \geq \frac{ 2 \nu \al^3}{1-\al^2} \ .
\]
Thus
\begin{equation}
\nu \leq \nu_0 (\al ) = \frac{\sqrt{32-3\al^6-17\al^4-16\al^2}}
{4(\al^2+1)(\al^2+4)} \quad (1/2 < \al < \al_0 )
\label{ext0}
\end{equation}
where $0.95 < \al_0 <1$ and $\al_0$ satisfies 
\[
32-3\al^6-17\al^4-16\al^2=0\ .
\]
For example, 
\[
\nu_0(0.5) \approx 0.244 \ , \quad \nu_0(0.95) \approx 0.0329 \ .
\]
For each fixed $\nu$, $\nu \leq \nu_0(\al )$, we have
\begin{equation}
f(0) > g(0)\ .
\label{prd4}
\end{equation}
From (\ref{prd3}) and (\ref{prd4}), we see that there is a $\la > 0$ 
such that (\ref{prd1}) is true. Next, we want to show that this eigenvalue 
is unique. Notice that
\begin{eqnarray}
(\la +\nu \al^2)^{-1} g(\la ) &=& \frac{2\al}{1-\al^2} \ , 
\label{prd5} \\
(\la +\nu \al^2)^{-1} f(\la ) &=& \frac{1}{(\la +\nu \al^2)a_{1} + 
\frac{1}{(\la +\nu \al^2)^{-1}a_{2}+ \cdots}} \ .
\label{prd6}
\end{eqnarray}
Since $(\la +\nu \al^2) a_{2n+1}$ ($n \geq 0$) is a strictly monotonically 
increasing function of $\la$ for $\la >0$, and $(\la +\nu \al^2)^{-1}a_{2n}$ 
($n \geq 1$) is a strictly monotonically decreasing function of $\la$ for 
$\la >0$, we see that $(\la +\nu \al^2)^{-1} f(\la )$ is a strictly 
monotonically decreasing function of $\la$ for $\la >0$. Thus the eigenvalue 
which satisfies
\[
(\la +\nu \al^2)^{-1} f(\la )=(\la +\nu \al^2)^{-1} g(\la )
\]
is unique. Similarly, $\nu \al^2 a_{2n+1}$ ($n \geq 0$) is a strictly 
monotonically increasing function of $\nu$ for $\nu >0$, and 
$(\nu \al^2)^{-1}a_{2n}$ ($n \geq 1$) is a constant function of $\nu$. 
Then $(\nu \al^2)^{-1}f(0)$ is a strictly monotonically decreasing function 
of $\nu$ for $\nu >0$. We know from above that when $\nu = \nu_0(\al )$,
\[
(\nu \al^2)^{-1}f(0) >(\nu \al^2)^{-1}g(0) = \frac{2\al}{1-\al^2}\ ,
\quad \mbox{constant in} \ \nu \ .
\]
From (\ref{prd2}), we have
\[
(\nu \al^2)^{-1}f(0)< \frac{1}{(\nu \al^2)a_1} = \frac{\al}{4(\al^2+1)^2}
\frac{1}{\nu^2}\ ,
\]
thus when $\nu >0$ is large enough
\[
(\nu \al^2)^{-1}f(0) <(\nu \al^2)^{-1}g(0)\ .
\]
Therefore, there is a unique $\nu_*(\al ) > \nu_0(\al )$, such that
\begin{equation}
(\nu_* \al^2)^{-1}f(0) = (\nu_* \al^2)^{-1}g(0)\ .
\label{ext1}
\end{equation}
Thus we have shown the following claims: When $\nu > \nu_*(\al)$, there 
is no eigenvalue of non-negative real part. 
When $\nu = \nu_*(\al )$, $\la =0$ is an eigenvalue, and all the rest 
eigenvalues have negative real parts. When $\nu < \nu_*(\al )$, there is 
a unique positive eigenvalue $\la (\nu )>0$, and all the rest 
eigenvalues have negative real parts. An estimate for $\nu_*$ can be obtained 
from (\ref{prd2}),
\[
\frac{1}{\frac{4(\al^2+1)^2}{\al} \nu_*^2 + \frac{\al^3 (\al^2+3)}
{4(\al^2+4)^2}} < \frac{2\al}{1-\al^2} < \frac{1}{\frac{4(\al^2+1)^2}{\al} 
\nu_*^2}\ .
\]
We have 
\[
\nu_0(\al ) < \nu_*(\al ) < \frac{1}{2(\al^2+1)}\sqrt{\frac{1-\al^2}{2}} \ ,
\]
where $\nu_0(\al )$ is given by (\ref{ext0}). For example,
\[
0.244029 <\nu_*(0.5) < 0.244949\ , \quad 
0.0329 < \nu_*(0.95) < 0.058\ .
\]
Next we want to show that the unique eigenvalue $\la (\nu )$ has the 
property that $\nu^{-1}\la (\nu )$ is a strictly monotonically decreasing 
function of $\nu$. Notice that
\begin{eqnarray*}
(\la +\nu \al^2) a_{2n+1} &=& \frac{4}{\al} \frac{\al^2+(2n+1)^2}
{\al^2+(2n+1)^2-1}\nu^2 \left (\frac{\la}{\nu}+\al^2 \right ) \\
& & \left (\frac{\la}{\nu}+[\al^2+(2n+1)^2]\right )\ , 
\quad (n\geq 0)\ , \\
(\la +\nu \al^2)^{-1} a_{2n} &=& \frac{4}{\al} \frac{\al^2+(2n)^2}
{\al^2+(2n)^2-1} \left [ 1 + \frac{(2n)^2}{\frac{\la}{\nu}+\al^2}
\right ]\ , \quad (n\geq 1)\ .
\end{eqnarray*}
Assume that $\nu^{-1}\la (\nu )$ is not a strictly monotonically 
decreasing function of $\nu$, then there is an interval in which 
$\nu^{-1}\la (\nu )$ is a strictly monotonically increasing or a constant 
function of $\nu$. In that interval, $[\la (\nu )+\nu \al^2]a_{2n+1}$ 
($n\geq 0$) is strictly monotonically increasing, and 
$[\la (\nu )+\nu \al^2]^{-1}a_{2n}$ ($n\geq 1$)
is monotonically decreasing. Thus $[\la (\nu )+\nu \al^2]^{-1}f(\la(\nu ))$ 
is a strictly monotonically decreasing function of $\nu$. On the other 
hand,
\[
[\la (\nu )+\nu \al^2]^{-1}g(\la(\nu ))= \frac{2\al}{1-\al^2}
\]
is a constant function of $\nu$. This contradiction shows that 
$\nu^{-1}\la (\nu )$ is a strictly monotonically decreasing function 
of $\nu$. Next we want to show that $\la (\nu )=O(1)$ as $\nu \ra 0^+$.
From (\ref{prd1}) and (\ref{prd2}), we have
\begin{eqnarray}
& & \frac{1}{\frac{4(\al^2+1)[\la +\nu (\al^2+1)]}{\al^3}+
\frac{\al (\al^2+3)}{4(\al^2+4)[\la +\nu (\al^2+4)]}} \non \\
&<& \frac{2\al}{1-\al^2}(\la + \nu \al^2) < 
\frac{\al^3}{4(\al^2+1)[\la +\nu (\al^2+1)]} \ .
\label{prd7}
\end{eqnarray}
Let $\nu \ra 0^+$, we have 
\begin{equation}
\sqrt{\frac{\al^2 (1-\al^2)}{8(\al^2+1)}-\frac{\al^4(\al^2+3)}
{16(\al^2+1)(\al^2+4)}} \leq \la \leq \sqrt{\frac{\al^2 (1-\al^2)}
{8(\al^2+1)}}\ .
\label{prd8}
\end{equation}
For the term under the square root on the left to be positive, we need
\begin{equation}
\al < \al_1 = \sqrt{\sqrt{\frac{59}{12}}-\frac{3}{2}} \approx 0.8469 \ .
\label{prd9}
\end{equation}
Thus when $\al \in (1/2,\al_1)$,
\[
\la (\nu )=O(1)\ , \quad \mbox{as}\ \nu \ra 0^+\ .
\]
Also from (\ref{prd7}), we have the estimate
\begin{eqnarray}
& & \sqrt{\frac{\al^2 (1-\al^2)}{8(\al^2+1)}-\frac{\al^4(\al^2+3)}
{16(\al^2+1)(\al^2+4)}}-\nu (\al^2+1) < \la (\nu ) \non \\
&<& \sqrt{\frac{\al^2 (1-\al^2)}{8(\al^2+1)}}- \nu \al^2 \ .
\label{prd95}
\end{eqnarray}
(b). $\nu = 0$. We have 
\[
a_n = \frac{4}{\al} \frac{\al^2+n^2}{\al^2+n^2-1}\la \ .
\]
As $|n| \ra \infty$, $a_n \ra \ta = \frac{4}{\al} \la$. As before, 
$\mbox{Re}\{ \la \} = 0$ and $|\la | \leq \frac{\al}{2}$ corresponds to 
a continuous spectrum. Next we show that outside the disc $|\la | 
\leq \frac{\al}{2}$, there is no eigenvalue. Outside the disc, $|\ta|>2$, 
\[
|a_0|>2\frac{\al^2}{1-\al^2} > \frac{2}{3}\ , \quad \mbox{for}\ 
\al \in (1/2,1)\ ,
\]
\[
|a_1|>2\left ( 1+ \frac{1}{\al^2}\right ) > 4\ , \quad \mbox{for}\ 
\al \in (1/2,1)\ ,
\]
and
\[
|a_n|>2\ , \quad \forall n \geq 2 \ .
\]
By a result of continued fraction \cite{LW92},
\[
\left |\frac{1}{a_{2} + \frac{1}{a_{3}+\cdots}}\right | \leq 1\ .
\]
Thus
\[
|g(\la )| > \frac{1}{3}\ ,
\]
while
\[
|f(\la )| \leq \frac{1}{|a_{1}| - \left |\frac{1}{a_{2} + 
\frac{1}{a_{3}+\cdots}}\right |} < \frac{1}{3}\ .
\]
This contradiction proves the claim. Thus the possible eigenvalues should 
lie in the region $\mbox{Re}\{ \la \} \neq 0$ and $|\la | \leq \frac{\al}{2}$.
Since the eigenvalues are symmetric with respect to the imaginary axis 
\cite{Li00}, we only need to consider the case $\mbox{Re}\{ \la \} > 0$. 
First we show that the possible eigenvalue must be real. Assume that 
$\mbox{Im}\{ \la \} > 0$, then
\[
\arg (-a_0) = \arg (a_1) = \arg (a_2) = \cdots > 0 \ .
\]
Thus 
\[
\left | \arg \left (\frac{1}{a_{2} + 
\frac{1}{a_{3}+\cdots}}\right ) \right | \leq \arg (a_1)\ .
\]
Then
\[
\arg \left ( a_1 + \frac{1}{a_{2} + 
\frac{1}{a_{3}+\cdots}}\right ) \neq -\arg (a_1)\ .
\]
Therefore
\[
\arg (f(\la )) \neq \arg (a_1) = \arg (-a_0) = \arg (g(\la ))\ .
\]
This is a contradiction. The case $\mbox{Im}\{ \la \} < 0$ is similar, and 
the claim is proved. When $\la >0$,
\[
\frac{1}{\la a_{1} + \frac{1}{\la^{-1}a_{2}}} <\la^{-1} f(\la ) < 
\frac{1}{\la a_{1}},
\]
i.e.
\begin{equation}
\frac{1}{\frac{4(\al^2+1)}{\al^3}\la^2 + \frac{\al (\al^2+3)}
{4(\al^2+4)}} <\la^{-1} f(\la ) < \frac{1}{\frac{4(\al^2+1)}{\al^3}\la^2}\ .
\label{prd10}
\end{equation}
Thus 
\[
\lim_{\la \ra +\infty} \la^{-1} f(\la ) = 0\ , \quad 
\lim_{\la \ra 0^+} \la^{-1} f(\la ) \geq \frac{4(\al^2+4)}{\al (\al^2+3)}\ .
\]
Notice that
\[
\la^{-1} g(\la ) =  \frac{ 2 \al }{1-\al^2} \ .
\]
We want to choose $\al$ such that
\[
\frac{ 2 \al }{1-\al^2} < \frac{4(\al^2+4)}{\al (\al^2+3)}\ .
\]
This condition is equivalent to the fact that the term under the square root 
on the left of (\ref{prd8}) being positive. Thus we have
\[
\al < \al_1 = \sqrt{\sqrt{\frac{59}{12}}-\frac{3}{2}} \approx 0.8469 \ ,
\]
which is the same with (\ref{prd9}). Thus when $\al \in (1/2,\al_1)$, there is 
a positive eigenvalue $\la_0$, 
\[
\la_0^{-1} f(\la_0 ) =\la_0^{-1} g(\la_0 ) = \frac{ 2 \al }{1-\al^2} \ .
\]
Since $\la a_{2n+1}$ ($n \geq 0$) is a strictly monotonically increasing 
function of $\la$ for $\la >0$, and $\la^{-1}a_{2n}$ ($n \geq 1$) is 
a constant function of $\la$, we see that $\la^{-1}f(\la )$ is a strictly 
monotonically decreasing function of $\la$ for $\la >0$. Thus the positive 
eigenvalue $\la_0$ is unique. From (\ref{prd10}), we have 
\[
\frac{1}{\frac{4(\al^2+1)}{\al^3}\la_0^2 + \frac{\al (\al^2+3)}
{4(\al^2+4)}} < \frac{ 2 \al }{1-\al^2} < \frac{1}{\frac{4(\al^2+1)}{\al^3}
\la_0^2}\ ,
\]
i.e.
\begin{equation}
\sqrt{\frac{\al^2 (1-\al^2)}{8(\al^2+1)}-\frac{\al^4(\al^2+3)}
{16(\al^2+1)(\al^2+4)}} < \la_0 < \sqrt{\frac{\al^2 (1-\al^2)}
{8(\al^2+1)}}\ ,
\label{prd11}
\end{equation}
which is similar to (\ref{prd8}). Next we want to show that 
\begin{equation}
\lim_{\nu \ra 0^+} \la (\nu ) = \la_0\ .
\label{prd12}
\end{equation}
Let $F(\nu , \la ) = f(\la )-g(\la )$, and let $F_N(\nu , \la )$ be the 
$N$-th truncation of $F(\nu , \la )$,
\[
F_N(\nu , \la )= \frac{a_0}{2} + \frac{1}{a_{1} + 
\frac{1}{a_{2}+\cdots +\frac{1}{a_n}}} \ .
\]
Notice that
\[
F(\nu , \la (\nu ) = 0\ , \quad \nu \in (0, \nu_*(\al ))\ , \quad 
\al \in (0.5, 0.8469)\ ; 
\]
\[
F(0, \la_0)= 0\ , \quad \al \in (0.5, 0.8469)\ ; 
\]
where $\la (\nu )$ and $\la_0$ have the estimates (\ref{prd95}) and 
(\ref{prd11}). Thus for any fixed $\al \in (0.5, 0.8469)$, 
($\nu,\la (\nu )$) lies in a compact set $[0,\dl]\times [c_1,c_2]$ where 
$\dl$, $c_1$ and $c_2$ are positive constants. Assume that (\ref{prd12}) 
is not true, then there is a sequence ($\nu_j,\la (\nu_j )$) such that 
\[
\lim_{\nu_j \ra 0^+} \la (\nu_j ) = \tla_0\ ,
\]
where 
\[
\tla_0 \neq \la_0 \quad \mbox{and} \quad \tla_0 \in [c_1,c_2]\ .
\]
By a result of continued fraction \cite{JT80},
\[
\lim_{N \ra \infty}F_N(\nu , \la )=F(\nu , \la )
\]
uniformly for $\nu \in [0,\dl]$ and $\la \in [c_1,c_2]$. Thus for any 
$\e >0$, there is an integer $N_0$ such that for any $N \geq N_0$,
\[
|F_N(\nu_j , \la (\nu_j )) -F(\nu_j , \la (\nu_j ))|<\e \ ,
\]
i.e.
\[
|F_N(\nu_j , \la (\nu_j ))|<\e \ , \quad \forall j \ .
\]
For any fixed $N$, let $j \ra +\infty$, we have
\[
|F_N(0 , \tla_0)|\leq \e \ .
\]
Let $N\ra +\infty$, we get
\[
|F(0 , \tla_0)|\leq \e \ .
\]
Since $\e$ is arbitrarily small, we have
\[
F(0 , \tla_0) = 0\ .
\]
By the uniqueness of the eigenvalue, $\tla_0 = \la_0$. This contradiction 
shows that (\ref{prd12}) is true.
\item When $\nu =0$, the union of all the continuous spatra 
\[
\left [ -i\frac{\al |\hk_1|}{2}, \ i\frac{\al |\hk_1|}{2} \right ]\ ,
\]
is $i\R$.
\end{enumerate}
\end{proof}

\subsection{Example 2}

Our second prime example is $\Om_* = \cos (x_1+x_2)$ which corresponds 
to $p=(1,1)^T$ and $\Ga = 1/2$ in (\ref{sm}). Here
\begin{eqnarray}
a_n &=& \frac{8}{(\hk_2-\hk_1)}\frac{(\hk_1+n)^2+(\hk_2+n)^2}
{(\hk_1+n)^2+(\hk_2+n)^2-2} \non  \\
& & \times \left \{ \la +\nu \left [ (\hk_1+n)^2+(\hk_2+n)^2 \right ]
\right \} \ .
\label{ape}
\end{eqnarray}
This example has two possible unstable eigenvalues. For shears, the number 
of unstable eigenvalues of 2D linear Euler is bounded by the number of 
lattices points inside the disc of radius $|p|$ \cite{LLS04}, which is 
even for square lattice. In \cite{Liu92} \cite{Liu95}, Vincent Liu studied 
another shear which has three possible unstable eigenvalues. Liu did a 
detailed calculation on the eigenvalues, which will be discussed in next 
subsection. Here we will do an even more detailed calculation. 
\begin{theorem}
The spectra of the 2D linear Euler operator $L$ have the following 
properties. 
\begin{enumerate}
\item The set $\{ \hk +np \}_{n \in \Z}$ has no intersection with the 
disc of radius $|p| = \sqrt{2}$. When $\nu >0$,  
there is no eigenvalue of non-negative real part. 
When $\nu = 0$, the entire spectrum is the continuous spectrum
\[
\left [ -i\frac{1}{4} |\hk_2 - \hk_1|, \ i\frac{1}{4} 
|\hk_2 - \hk_1|\right ]\ .
\]
\item $\hk = p =(1,1)^T$. The spectrum consists of the eigenvalues 
\[
\la = - \nu 2n^2 \ , \quad n \in \Zo \ .
\]
The eigenfunctions are the Fourier modes
\[
\tom_{np} e^{i n p\cdot x} + \ \mbox{c.c.}\ \ , \quad \forall \tom_{np} \in 
\C\ , \quad n \in \Zo \ .
\]
As $\nu \ra 0^+$, the eigenvalues are dense on the negative half of the real 
axis $(-\infty, 0]$. Setting $\nu =0$, the only eigenvalue is $\la = 0$ of 
infinite multiplicity with the same eigenfunctions as above.
\item $\hk =(-1,1)^T$. When $\nu >0$,  
there is no eigenvalue of non-negative real part. 
When $\nu = 0$, the entire spectrum is the continuous spectrum
$[-i\frac{1}{2}, i\frac{1}{2}]$. A special eigenvalue is $\la = -2\nu$ ( 
when $\nu = 0$, this eigenvalue $\la =0$ is embedded in the continuous 
spectrum). 
\item $\hk =(0,1)^T$. When $\nu >0$, in the half plane 
$\mbox{Re}\{ \la \} \geq -\nu$, there is a unique pair of eigenvalues 
$\la$ and $\bar{\la}$ such that
\[
-\nu < \ \mbox{Re}\{ \la \} \ < \frac{1}{4}\sqrt{\frac{3}{20}+(8\nu )^2}
-2 \nu \ , 
\]
\[
\frac{1}{8}\left ( 1-\sqrt{\frac{3}{5}}\right ) < \ \mbox{Im}\{ \la \} \ <
\frac{1}{8}\left ( 1+\sqrt{\frac{3}{5}}\right )\ .
\]
When $\nu = 0$, $[-i\frac{1}{4}, i\frac{1}{4}]$ is a continuous spectrum. 
If there is an eigenvalue of positive real part, then there is a quadruplet
($\la$, $\bar{\la}$ $-\la$ $-\bar{\la}$) where
\[
0< \ \mbox{Re}\{ \la \} \ < \frac{1}{16}\sqrt{\frac{3}{5}}\ , 
\]
\[
\frac{1}{8}\left ( 1-\sqrt{\frac{3}{5}}\right ) < \ \mbox{Im}\{ \la \} \ <
\frac{1}{8}\left ( 1+\sqrt{\frac{3}{5}}\right )\ .
\]
\item Finally, when $\nu = 0$, the union of all the above pieces of 
continuous spectra is the imaginary axis $i\R$.
\end{enumerate}
\label{PTHM}
\end{theorem}
\begin{remark}
For Case 4, numerical computation indicates that when $\nu > 0$, the 
real part of the eigenvalue $\mbox{Re}\{ \la (\nu ) \} \geq c >0$, where 
$c$ is independent of $\nu$. When $\nu = 0$, numerical computation 
indicates that the eigenvalue $\la (0)$ indeed exists and $\mbox{Re}\{ 
\la (0) \} \geq c >0$ (the same $c$ as above). Numerical computation 
also indicates that as $\nu \ra 0^+$, $\la (\nu ) \ra \la (0)$ 
which can be proved given the above facts. Cases 2 and 3 indicate that 
$\mbox{Re}\{ \la^s_1 \}$ in Assumption 1 is at least $O(\nu )$ as 
$\nu \ra 0^+$.
\end{remark}
\begin{proof}
We will give the proof case by case.
\begin{enumerate}
\item The case (a). $\nu > 0$. If $\mbox{Re}\{ \la \} \geq 0$, then 
all the $\mbox{Re}\{ a_n \}$'s are nonzero and have the same sign.
The real parts of the right and left hand sides of (\ref{mch}) are 
nonzero but of different signs. Thus there is no eigenvalue of 
non-negative real part. 
The case (b). $\nu = 0$. If $\mbox{Re}\{ \la \} \neq 0$, then all 
the $\mbox{Re}\{ a_n \}$'s are nonzero and have the same sign. Again 
(\ref{mch}) can not be satisfied. If $\mbox{Re}\{ \la \} = 0$, let 
$\ta = \lim_{n \ra \infty} a_n = \frac{8}{\hk_2- \hk_1} \la$, then
\[
\ta z_n = \frac {|\hk +np|^2 -2}{|\hk +np|^2}(z_{n+1}-z_{n-1})\ ,
\]
where $0< |\hk +np|^2 -2 < |\hk +np|^2$, $\forall n \in \Z$. By using 
$\ell_2$ norm of $\{ z_n \}_{n \in \Z}$, one sees that possible 
eigenvalues have to satisfy $|\ta| \leq 2$. As mentioned before, 
$\mbox{Re}\{ \la \} = 0$ and $|\ta | \leq 2$ correspond to a continuous 
spectrum \cite{Li00}, which is the interval
\begin{equation}
\left [ -i\frac{1}{4}|\hk_2- \hk_1|, \ i\frac{1}{4}|\hk_2- \hk_1| \right ]\ .
\label{piece}
\end{equation}
\item In this case, one has
\[
[\la +\nu 2(n+1)^2]\tom_{(n+1)p} = 0 \ .
\]
The claims follow immediately.
\item In this case, denote $\tom_{\hk +np}$ simply by $\om_n$, one has
\[
[\la +2\nu (n^2+1)]\om_n = \frac{1}{4} \frac{(n+1)^2}{(n+1)^2+1} 
\om_{n+1} - \frac{1}{4} \frac{(n-1)^2}{(n-1)^2+1} 
\om_{n-1} \ .
\]
This system decouples at $n=0$. A few equations around $n=0$ are 
\begin{eqnarray}
(\la +10 \nu ) \om_{-2} &=& \frac{1}{8}\om_{-1}-\frac{9}{40}\om_{-3} \ ,
\non \\
(\la +4 \nu ) \om_{-1} &=& -\frac{1}{5}\om_{-2}\ , \label{deg1} \\
(\la +2 \nu ) \om_{0} &=& \frac{1}{8}\om_{1}-\frac{1}{8}\om_{-1} \ ,
\non \\
(\la +4 \nu ) \om_{1} &=& \frac{1}{5}\om_{2}\ , \label{deg2} \\
(\la +10 \nu ) \om_{2} &=&\frac{9}{40}\om_{3}-\frac{1}{8}\om_{1} \ . \non
\end{eqnarray}
Notice that
\[
a_n = 4 \left ( 1+\frac{1}{n^2}\right ) [\la +2\nu (n^2+1)]\ .
\]
From (\ref{deg2}), one gets
\begin{equation}
w_2 = z_2/z_1 =8 (\la +4 \nu )\ .
\label{deg3}
\end{equation}
From (\ref{CF1})and (\ref{deg3}), one has 
\begin{equation}
8(\la +4 \nu ) = - \frac{1}{a_2+\frac{1}{a_3+\cdots}}\ . 
\label{deg4}
\end{equation}
(a). $\nu >0$. If $\mbox{Re}\{ \la \} \geq 0$, then $\mbox{Re}\{ a_n \}
>0$ for $n \geq 2$. Thus the real part of the right hand side of (\ref{deg4})
is negative, while the real part of its left hand side is positive. Thus 
there is no eigenvalue of non-negative real part. (b). $\nu =0$. 
If $\mbox{Re}\{ \la \} \geq 0$, then $\mbox{Re}\{ a_n \}$ ($n\geq 2$) has 
the same fixed sign with $\la$, again (\ref{deg4}) leads to a contradiction.
If $\mbox{Re}\{ \la \} = 0$, let $\ta = \lim_{n\ra \infty} a_n = 4\la$, 
then 
\begin{eqnarray*}
\ta z_1 &=& \frac{1}{2} z_2 \ , \\
\ta z_n &=& \frac{n^2}{n^2+1} (z_{n+1}-z_{n-1})\ , \quad n \geq 2 \ .
\end{eqnarray*}
By using 
$\ell_2$ norm of $\{ z_n \}_{n \geq 1}$, one sees that possible 
eigenvalues have to satisfy $|\ta| \leq 2$. As mentioned before, 
$\mbox{Re}\{ \la \} = 0$ and $|\ta | \leq 2$ correspond to a continuous 
spectrum \cite{Li00}, which is the interval $[-i\frac{1}{2}, i\frac{1}{2}]$.
Similarly from (\ref{deg1}), one gets
\begin{equation}
w_{-1} =z_{-1}/z_{-2}=-\frac{1}{8(\la +4 \nu )}\ .
\label{deg5}
\end{equation}
From (\ref{CF1})and (\ref{deg5}), one has 
\[
-\frac{1}{8(\la +4 \nu )} = a_{-2}+ \frac{1}{a_{-3}+\frac{1}
{a_{-4}+\cdots}}\ . 
\]
Notice that $a_{-n}=a_n$, one gets the same conclusion as above. Finally, 
by choosing $\om_n = 0$ if $n\neq 0$, one gets the eigenvalue $\la = -2\nu$ 
with the eigenfunction 
\[
\tom_{\hk} e^{i\hk \cdot x} + \ \mbox{c.c.}\ \ , \quad \forall 
\tom_{\hk} \in \C \ .
\]
\item In this case,
\[
a_n = 8 \frac{n^2 +(n+1)^2}{n^2 +(n+1)^2-2} \left \{ \la +\nu \left [ 
n^2 +(n+1)^2\right ] \right \} \ .
\]
Thus $a_{-(n+1)} = a_n$. Equation (\ref{mch}) is reduced to 
\begin{equation}
a_0 + \frac{1}{a_1+\frac{1}{a_2+\cdots}} = \pm i \ .
\label{rde}
\end{equation}
The first a few $a_n$'s are 
\begin{eqnarray*}
a_0 &=& 8\cdot (-1) \cdot (\la +\nu )\ , \\
a_1 &=& 8\cdot (\frac{5}{3}) \cdot (\la +5\nu )\ , \\
a_2 &=& 8\cdot (\frac{13}{11}) \cdot (\la +13\nu ) \ , \\
a_3 &=& 8\cdot (\frac{25}{23}) \cdot (\la +25\nu ) \ .
\end{eqnarray*}
(a). $\nu >0$. Consider the region $\mbox{Re}\{ \la \} \geq -\nu$, in which 
\[
|a_n| \geq 8|\la +\nu |\ , \quad \forall n \geq 0 \ .
\]
By using 
$\ell_2$ norm of $\{ z_n \}_{n \in \Z}$ in (\ref{PP}), one sees that possible 
eigenvalues in this region have to satisfy
\[
|\la +\nu |\leq \frac{1}{4}\ .
\]
Next in the region 
\[
\D = \bigg \{ \la \ \bigg | \ \mbox{Re}\{ \la \} \geq -\nu, \ 
|\la +\nu |\leq \frac{1}{4} \bigg \} \ ,
\]
we will use Rouch\'e's theorem to track the eigenvalues. Let
\begin{eqnarray}
f(\la ) &=& a_0 + \frac{1}{a_1+\frac{1}{a_2+\cdots}} + i \ ,
\label{rtf} \\
g(\la ) &=& a_0 + \frac{1}{a_1}+ i \ . \non
\end{eqnarray}
To apply the Rouch\'e's theorem, one needs to show that
\begin{equation}
|f(\la )-g(\la )| <|f(\la )| +|g(\la )|\ , \quad \la \in \pa \D \ .
\label{Rou1}
\end{equation}
If this is not true, then there is a $\la \in \pa \D$ and a $\dl$ 
($0\leq \dl \leq \infty$) such that
\[
f(\la )= -\dl g(\la )\ ,
\]
i.e.
\begin{equation}
a_0 + i + \frac{\frac{1}{1+\dl}}{a_1+\frac{1}{a_2+\cdots}}
+\frac{\dl}{1+\dl}\frac{1}{a_1} = 0 \ .
\label{Rou2}
\end{equation}
On the part of the boundary $\pa \D$: $\mbox{Re}\{ \la \} = -\nu$,
\[
\mbox{Re}\{ a_0 \} = 0 \ , \quad \mbox{Re}\{ a_n\} > 0 \ , \quad n \geq 1 \ .
\]
By taking the real part of (\ref{Rou2}), one sees that (\ref{Rou2}) 
can not be satisfied. On the other part the boundary $\pa \D$: 
$|\la +\nu | =\frac{1}{4}$,
\begin{eqnarray*}
& & |a_0|=2\ , \quad |a_1| \geq \frac{10}{3}\ , \quad  |a_2| 
\geq \frac{26}{11}\ , \\
& & |a_n| \geq 2\ , \quad \forall n \geq 3\ .
\end{eqnarray*}
By a result of continued fraction \cite{LW92},
\[
\left |\frac{1}{a_{3} + \frac{1}{a_{4}+\cdots}}\right | \leq 1\ .
\]
Thus
\begin{eqnarray*}
\left | \frac{\frac{1}{1+\dl}}{a_1+\frac{1}{a_2+\cdots}} \right |
&\leq& \frac{1}{|a_1|-\left |\frac{1}{a_2+\frac{1}{a_3+\cdots}}\right |} \\
&\leq& \frac{1}{|a_1|-\frac{1}{|a_2|-
\left |\frac{1}{a_3+\frac{1}{a_4+\cdots}}\right |}} \\
&\leq& \frac{1}{\frac{10}{3}-\frac{1}{\frac{26}{11}-1}} = \frac{15}{39}\ .
\end{eqnarray*}
The rest of (\ref{Rou2}) has the estimate:
\begin{eqnarray*}
& & \left | a_0 + i + \frac{\dl}{1+\dl}\frac{1}{a_1} \right | 
\geq |a_0| - 1 - \frac{\dl}{1+\dl}\frac{1}{|a_1|} \\
& & \geq 2 - 1 - \frac{\dl}{1+\dl} \frac{3}{10} \geq \frac{7}{10} \ .
\end{eqnarray*}
Therefore, (\ref{Rou2}) can not be satisfied. By Rouch\'e's theorem, 
inside $\D$, $f(\la )$ and $g(\la )$ have the same number of roots.
$g(\la )$ has two roots
\[
\la =-3\nu + \frac{1}{16}i \pm \frac{1}{16}\sqrt{\frac{7}{5} + 16(8\nu )^2
+i 64 \nu }\ ,
\]
and one of which is in $\D$ (at least when $\nu$ is small enough). Thus 
$f(\la )$ has one root in $\D$. From (\ref{rtf}), this root satisfies
\begin{equation}
\mbox{Re} \left \{ a_1 +\frac{1}{a_0+i}\right \} < 0 \ ,
\label{pst}
\end{equation}
which leads to 
\[
-\nu < \ \mbox{Re}\{ \la \} \ < \frac{1}{4}\sqrt{\frac{3}{20}+(8\nu )^2}
-2 \nu \ , 
\]
\[
\frac{1}{8}\left ( 1-\sqrt{\frac{3}{5}}\right ) < \ \mbox{Im}\{ \la \} \ <
\frac{1}{8}\left ( 1+\sqrt{\frac{3}{5}}\right )\ .
\]
(b). $\nu =0$. Then
\[
|a_n| \geq 8|\la |\ , \quad \forall n \geq 0 \ .
\]
Again by using $\ell_2$ norm of $\{ z_n \}_{n \in \Z}$ in (\ref{PP}), one 
sees that possible eigenvalues in this region have to satisfy
\[
|\la |\leq \frac{1}{4}\ .
\]
In fact, if there is an eigenvalue with positive real part, then 
(\ref{pst}) is true, which leads to 
\[
0< \ \mbox{Re}\{ \la \} \ < \frac{1}{16}\sqrt{\frac{3}{5}}\ , 
\]
\[
\frac{1}{8}\left ( 1-\sqrt{\frac{3}{5}}\right ) < \ \mbox{Im}\{ \la \} \ <
\frac{1}{8}\left ( 1+\sqrt{\frac{3}{5}}\right )\ .
\]
As proved in \cite{Li00}, such eigenvalues (if exist) come in quadruplet
($\la$, $\bar{\la}$ $-\la$ $-\bar{\la}$). As in Case 3 above, 
$[-i\frac{1}{4}, i\frac{1}{4}]$ is a continuous spectrum \cite{Li00}.
\item When $\nu =0$, the union of all the continuous spectra
\[
\left [ -i\frac{1}{4}|\hk_2-\hk_1|, \ i\frac{1}{4}|\hk_2-\hk_1|\right ]\ ,
\]
is $i\R$.
\end{enumerate}
\end{proof}

\subsection{Example 3}

An example studied in details by Vincent Liu \cite{Liu92} \cite{Liu95}
is $\Om_* = -\frac{\sqrt{2}}{\pi}\cos (2x_2)$ which corresponds 
to $p=(0,2)^T$ and $\Ga = -\frac{1}{\sqrt{2}\pi}$ in (\ref{sm}). Here
\begin{eqnarray*}
a_n &=& -4\sqrt{2}\pi \frac{\hk_1^2+(\hk_2+2n)^2}
{\hk_1 \left [\hk_1^2+(\hk_2+2n)^2-4\right ] } \\
& &\times \left \{ \la +\nu \left [  \hk_1^2+(\hk_2+2n)^2\right ]\right \} \ .
\end{eqnarray*}
\begin{theorem}
The spectra of the 2D linear Euler operator $L$ have the following 
properties. 
\begin{enumerate}
\item The set $\{ \hk +np \}_{n \in \Z}$ has no intersection with the 
disc of radius $|p| = 2$. When $\nu >0$, 
there is no eigenvalue of non-negative real part. 
When $\nu = 0$, the entire spectrum is the continuous spectrum
\[
\left [ -i\frac{|\hk_1|}{2\sqrt{2}\pi}, \ i\frac{|\hk_1|}{2\sqrt{2}\pi}
\right ]\ .
\]
\item $\hk =(0,1)^T$ or $(0,2)^T$. The spectrum consists of the eigenvalues 
\[
\la = - \nu [  \hk_1^2+(\hk_2+2n)^2 ] \ , \quad n \in \Z \ ,
\]
where for $\hk =(0,2)^T$, $n \neq -1$. 
The eigenfunctions are the Fourier modes
\[
\tom_{\hk +np} e^{i(\hk +np)\cdot x} + \ \mbox{c.c.}\ \ , 
\quad \forall \tom_{\hk +np} \in 
\C\ , \quad n \in \Z \ .
\]
As $\nu \ra 0^+$, the eigenvalues are dense on the negative half of the real 
axis $(-\infty, 0]$. Setting $\nu =0$, the only eigenvalue is $\la = 0$ of 
infinite multiplicity with the same eigenfunctions as above.
\item $\hk =(2,0)^T$. When $\nu >0$, 
there is no eigenvalue of non-negative real part. 
When $\nu = 0$, the entire spectrum is the continuous spectrum
$[-i\frac{1}{\sqrt{2}\pi}, i\frac{1}{\sqrt{2}\pi}]$. A special eigenvalue 
is $\la = -4\nu$ (when $\nu = 0$, this eigenvalue $\la =0$ is embedded in 
the continuous spectrum). 
\item $\hk = (1,0)^T$. There is a unique $\nu_*$,
\[
\frac{1}{10\pi }\sqrt{\frac{89}{34}} < \nu_* < \frac{1}{10\pi }\sqrt{3}\ .
\]
When $\nu > \nu_*$, there is no eigenvalue of non-negative real part. 
When $\nu = \nu_*$, $\la =0$ is an eigenvalue, and all the rest 
eigenvalues have negative real parts. When $\nu < \nu_*$, there is 
a unique positive eigenvalue $\la (\nu )>0$, and all the rest 
eigenvalues have negative real parts. $\nu^{-1} \la (\nu )$ is a strictly 
monotonically decreasing function of $\nu$. $\la (\nu )$ has the estimate
\[
\frac{1}{\pi}\sqrt{\frac{89}{680}} - 5 \nu < \la (\nu ) 
< \frac{1}{\pi}\sqrt{\frac{3}{20}} - \nu \ .
\]
In particular, as $\nu \ra 0^+$, $\la (\nu ) =O(1)$.
When $\nu =0$, we have only two eigenvalues
$\la_0$ and $-\la_0$, where $\la_0$ is positive
\[
\frac{1}{\pi}\sqrt{\frac{89}{680}}  < \la_0 
< \frac{1}{\pi}\sqrt{\frac{3}{20}} \ .
\]
The rest of the spectrum is a continuous spectrum 
$[-i\frac{1}{2\sqrt{2}\pi}, i\frac{1}{2\sqrt{2}\pi}]$. Moreover,
\[
\lim_{\nu \ra 0^+} \la (\nu ) = \la_0 \ .
\]
\item $\hk =(1,1)^T$. When $\nu >0$, in the half plane 
$\mbox{Re}\{ \la \} \geq -2\nu$, there is a unique pair of eigenvalues 
$\la$ and $\bar{\la}$ such that
\[
-2\nu < \ \mbox{Re}\{ \la \} \ <  \frac{1}{\pi}\sqrt{\frac{3}{40}}\ , 
\]
\[
\frac{\sqrt{5} -\sqrt{3}}{2\sqrt{10}\pi} < \ \mbox{Im}\{ \la \} \ <
\frac{\sqrt{5} +\sqrt{3}}{2\sqrt{10}\pi} \ .
\]
When $\nu = 0$, $[-i\frac{1}{2\sqrt{2}\pi}, i\frac{1}{2\sqrt{2}\pi}]$ is 
a continuous spectrum. 
If there is an eigenvalue of positive real part, then there is a quadruplet
($\la$, $\bar{\la}$ $-\la$ $-\bar{\la}$) where
\[
0< \ \mbox{Re}\{ \la \} \ < \frac{1}{\pi}\sqrt{\frac{3}{40}}\ , 
\]
\[
\frac{\sqrt{5} -\sqrt{3}}{2\sqrt{10}\pi} < \ \mbox{Im}\{ \la \} \ <
\frac{\sqrt{5} +\sqrt{3}}{2\sqrt{10}\pi} \ .
\]
\item Finally, when $\nu = 0$, the union of all the above pieces of 
continuous spectra is the imaginary axis $i\R$.
\end{enumerate}
\end{theorem}
\begin{remark}
Except for Case 4, the rest of the proof is similar to that of Theorem 
\ref{PTHM}. Overall, Liu \cite{Liu92} \cite{Liu95} did not realize the 
continuous spectrum. For Case 4, Liu \cite{Liu92} \cite{Liu95} gave an 
elegant proof. For Case 5, Liu \cite{Liu92} \cite{Liu95} claimed more 
than what were actually proved. In \cite{Liu95}, the sketch after line 9, 
pp.472 can be realized as shown in the proof of Case 4 of Theorem 
\ref{PTHM} above. But Equation (39) on the same page does not imply 
$\mbox{Re} \{\eta_1(\nu )\} >0$ in Equation (5) on pp.467. Arguments between 
lines 5 and 8 on pp.483 are not solid, therefore, the existence of $\eta_1$ 
in Part B of Theorem 2 on pp.467 was not proved. The sketch between lines 
22 and 25 on pp.484 is not completed, therefore, (14) in Theorem 3 on 
pp.468 was not proved. Proving these claims seems tricky.
\end{remark}
\begin{remark}
Cases 2 and 3 indicates that $\mbox{Re} \{ \la_1^s \}$ in Assumption 1 is 
at least $O(\nu )$ as $\nu \ra 0^+$. Case 4 shows the existence of an unstable
eigenvalue. In fact, by the property of continued fraction, the corresponding 
eigenfunction $\Om^u(\nu )$ also converges as $\nu \ra 0^+$. Thus 
$\| \Om^u \|_{\ell +1}/\| \Om^u \|_{\ell}$ has a uniform bound 
as $\nu \ra 0^+$.
\end{remark}
\begin{remark}
In an effort to minimize the number of unstable modes to one and to have 
Case 4 not Case 5 in the above theorem, we need to study the rectangular 
periodic domain $[0,2\pi /\al ] \times [0, 2\pi ]$ where $1/2 < \al <1$, and 
the steady state $\Om_* = \cos x_2$. This is our Example 1 studied before.
\end{remark}

\end{document}